\documentclass[%
 reprint,
 showpacs,preprintnumbers,
 amsmath,
 amssymb,
 aps,
 pre,
]{revtex4-1}

\usepackage{graphicx}
\usepackage{dcolumn}
\usepackage{bm}
\usepackage{hyperref}
\usepackage{accents}
\usepackage{mathrsfs}
\usepackage{xcolor}
\usepackage[normalem]{ulem}
\usepackage{amsthm}

\newcommand{\ab}{{\alpha,\beta}}

\begin{document}

\title{Identifying Invariant Ergodic Subsets and Barriers to Mixing by Cutting and Shuffling:\texorpdfstring{\\}{\space}  Study in a Bi-rotated Hemisphere}

\author{Thomas F. Lynn}
\affiliation{Department of Engineering Sciences and Applied Mathematics, Northwestern University, Evanston, Illinois 60208, USA}

\author{Julio M. Ottino}
\affiliation{Department of Chemical and Biological Engineering, Department of Mechanical Engineering, Northwestern University, Evanston, Illinois 60208, USA}
\affiliation{The Northwestern Institute on Complex Systems (NICO), Northwestern University, Evanston, Illinois 60208, USA}

\author{Paul B. Umbanhowar}
\affiliation{Department of Mechanical Engineering, Northwestern University, Evanston, Illinois 60208, USA}

\author{Richard M. Lueptow}
\email{r-lueptow@northwestern.edu}
\affiliation{Department of Mechanical Engineering, Department of Chemical and Biological Engineering, Northwestern University, Evanston, Illinois 60208, USA}
\affiliation{The Northwestern Institute on Complex Systems (NICO), Northwestern University, Evanston, Illinois 60208, USA}%

\date{\today}

\begin{abstract}
Mixing by cutting-and-shuffling can be mathematically described by the dynamics of piecewise isometries (PWIs), higher dimensional analogs of one-dimensional interval exchange transformations. In a two-dimensional domain under a PWI, the exceptional set, \(\bar{E}\), which is created by the accumulation of cutting lines (the union of all iterates of cutting lines and all points that pass arbitrarily close to a cutting line), defines where mixing is possible but not guaranteed.  There is structure within \(\bar{E}\) that directly influences the mixing potential of the PWI. Here we provide new computational and analytical formalisms for examining this structure by way of measuring the density and connectivity of \(\varepsilon\)-fattened cutting lines that form an approximation of \(\bar{E}\). For the example of a PWI on a hemispherical shell studied here, this approach reveals the subtle mixing behaviors and barriers to mixing formed by invariant ergodic subsets (confined orbits) within the fractal structure of the exceptional set. Some PWIs on the shell have provably non-ergodic exceptional sets, which prevent mixing, while others have potentially ergodic exceptional sets where mixing is possible since ergodic exceptional sets have uniform cutting line density.  For these latter exceptional sets, we show the connectivity of orbits in the PWI map through direct examination of orbit position and shape and through a two-dimensional return plot to explain the necessity of orbit connectivity for mixing.
\end{abstract}

\maketitle


\section{Introduction \label{sec:introduction}}

The mathematical foundation of cut-and-shuffle mixing is the piecewise isometry (PWI), which cuts a domain into pieces that are rearranged to reform the original domain \cite{Ashwin1997a, Ashwin1997, Goetz2000, Adler2001, Goetz2004, Smith2019}. Mixing by cutting-and-shuffling arises in several natural systems such as granular materials \cite{Sturman2012, Juarez2010, Juarez2012, Park2017, Smith2017, Umbanhowar2013}, valved fluid flows \cite{Jones1988, Smith2016a},  and imbricate thrust faults in geology \cite{Boyer1982, Bell1983, Butler1982}. In these systems and others, cutting-and-shuffling is not the only mixing mechanism (chaotic advection and diffusion can also be present), but it is presently the least understood. Here we present techniques to investigate the mixing of a particular class of hemispherical PWIs related to the mixing of granular particles in a spherical tumbler \cite{Park2016, Park2017, Juarez2010, Juarez2012}. Our goal is not to study the wide range of PWIs, which can occur in one, two, and three dimensions (or higher) \cite{Smith2019}, but instead to develop computational and analytical approaches to identify invariant ergodic subsets that result in barriers to mixing. And, more importantly, the results we will describe are largely independent of the choice of (almost everywhere invertible) PWI.

PWIs act isometrically (i.e., as a distance preserving transformation) on each of a finite number of pieces, or ``atoms,'' \(P_i\) of a domain. For this paper, PWI dynamics are examined near the boundaries of these pieces where the map is discontinuous. Following previous work \cite{Kahng2009, Park2017, Lynn2019, Smith2017}, the PWI here is defined as a multi\-valued map acting on the power set (the set of subsets, which allows multi\-valued boundaries) such that the different, but overlapping, boundaries of the pieces of the domain are allowed to have different actions under the map. This is to say that, under the map, atoms that touch each other retain a copy of their mutual boundary under the map. This is an identical treatment to our previous papers \cite{Park2017, Lynn2019}, although it was not emphasized in the arguments presented there. Atom boundaries are essential elements of some of the arguments presented in this paper, and it is a point that requires special care. Other treatments define dynamics on cutting lines such that the map is not multi\-valued \cite{Ashwin1997, Goetz2000, Goetz2001, Scott2001, Scott2003, Ashwin2002, Goetz2003, Goetz2004}, and single-valued cutting lines are sufficient for some arguments here (the computational methods are ultimately indifferent to the treatment of the boundaries).

\begin{figure}
	\includegraphics[width=0.45\textwidth]{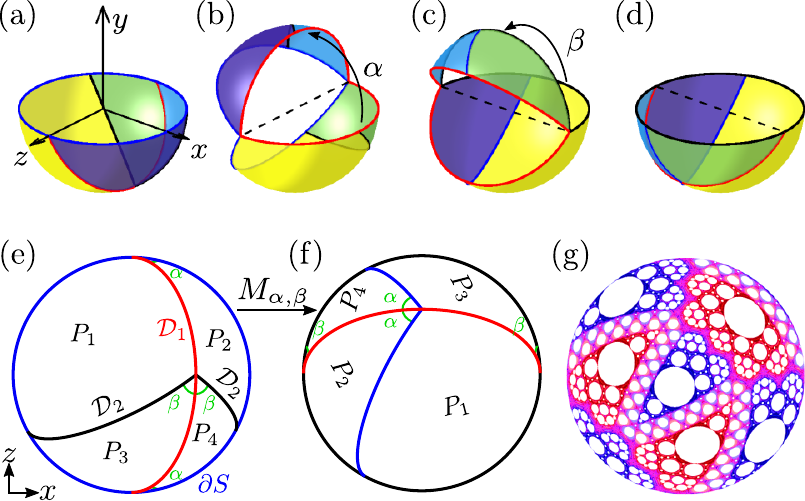}
  \caption{(a-d) Illustration of the rotations that define the hemispherical PWI studied here. (a) Initial condition of the hemispherical shell (\(S\)). (b) Rotation about the \(z\)-axis by \(\alpha\) (arrow shows direction only).
  (c) Rotation about the \(x\)-axis by \(\beta\).
  (d) Re-formed \(S\).
  (e,f) Equivalent PWI for the \((\alpha = 57^{\circ}, \beta = 57^{\circ})\) protocol results in four atoms, \(P_1, P_2, P_3, P_4 \), with cutting lines \(\mathcal{D}_1\) and \(\mathcal{D}_2\).
  (g) Accumulation of pre- and post-images of \(\varepsilon\)-fattened (\(\varepsilon = 10^{-3}\)) cutting lines \(\mathcal{D}\) after \(N = 20\,000\) iterations, viewed from the negative \(y\)-axis, approximating the true exceptional set \(\bar{E}\) as \(\tilde{E}_{\varepsilon,N}\) using color to indicate cutting line density as explained in Sec.~\ref{sec:natural_invariant}. A Lambert azimuthal equal area projection is used \cite{Snyder1987}.
  }
  \label{fig:pwi_demo}
\end{figure}

The specific system we study is inspired by a physically realizable spherical tumbler half-filled with granular matter and rotated in a periodic sequence about two perpendicular axes \cite{Zaman2018,Yu2019,Juarez2010}, though the approach described here to identify barriers to mixing could be applied to any PWI. In the equivalent hemispherical shell PWI example we consider here \cite{Juarez2012,Park2017,Smith2017,Park2016,Lynn2019,Scott2001,Scott2003,Zaman2018,Yu2019,Park2016a}, the lower unit hemispherical shell, \(S = \{\bm{x} = (x,y,z) \in \mathbb{R}^3 : \|\bm{x}\| = 1, y\le 0\}\), in Fig.~\ref{fig:pwi_demo}(a) is rotated by angle \(\alpha \in [0,\pi)\) about the \(z\)-axis and then by angle \(\beta \in [0,\pi)\) about the \(x\)-axis, as shown in Fig.~\ref{fig:pwi_demo}(b,c) to remake \(S\) in Fig.~\ref{fig:pwi_demo}(d). The equator provides the cutting action as the pieces of the hemisphere that are flipped across the boundary separate from neighboring points. This rotational procedure, i.e.\ Fig.~\ref{fig:pwi_demo}(a-d), motivates the equivalent PWI map \(M_\ab\), referred to as \(M\) except where \(\ab\) need to be specified, corresponding to the \textit{protocol} \((\alpha,\beta)\) shown in Fig.~\ref{fig:pwi_demo}(e,f).

In this PWI, the hemispherical domain \(S\) is split into at most four atoms, \(\{P_1,P_2,P_3,P_4\}\), as shown in Fig.~\ref{fig:pwi_demo}(e), which are mutually disjoint but overlap at their shared boundaries, \(\mathcal{D}_1\) and \(\mathcal{D}_2\). These atoms are rearranged according to the rotations in Fig.~\ref{fig:pwi_demo}(a-d) as shown in Fig.~\ref{fig:pwi_demo}(f). Repeating the PWI, i.e.\ cutting-and-shuffling over and over, reveals the \textit{singular set}, \(E\), which is defined as \(E = \bigcup_{n=0}^\infty M^n(\mathcal{D}_1 \cup \mathcal{D}_2)\) and is the accumulation of cuts produced by the PWI map \(M\) (all multi\-valued points of the map and the union of forward iterates of cutting lines) \cite{Kahng2009}. As previously defined, the singular set and its limit points define the closure of \(E\) called the \textit{exceptional set}, \(\bar{E}\), which is \(E\) and all points that pass arbitrarily close to \(E\) under the PWI map \footnote{The limit points of a set define all the points for which there is a sequence of points within the set that approach the limit point with arbitrary closeness. In this case, \(E\) is the set of all points \(\bm{x}\) for which the minimum distance between \(\{M^n(\bm{x}): -\infty < n < \infty\}\), i.e.\ the \textit{orbit} of \(\bm{x}\), and \(\mathcal{D}\) is exactly zero; \(\min_n d(M^n(\bm{x},\mathcal{D}))=0\) for a distance metric \(d\). The limit points of \(E\) are all points for which the infimum of this distance is zero, but this distance is never exactly zero; \(\inf_n d(M^n(\bm{x},\mathcal{D}))=0\) and \(d(M^n(\bm{x},\mathcal{D})) > 0\).}.

For generic protocols, it has been conjectured that \(\bar{E}\) is a fat fractal \cite{Farmer1983, Umberger1985}, and there is strong numerical evidence that \(\bar{E}\) has non-zero measure for almost all protocols \cite{Ashwin1997a, Ashwin1997, Bruin2003, Scott2003, Park2017} (that is, a finite fraction of the domain passes arbitrarily close to the cutting lines under the PWI map) except for those protocols that produce polygonal tilings \cite{Lynn2019, Smith2017}. Instead of an exact representation of \(\bar{E}\),  an \(N\)-iterate approximation \cite{Lynn2019} to the exceptional set using \(\varepsilon\)-fattened cutting lines, \(\tilde{E}_{\varepsilon,N}\), is used here; for the PWI in Fig.~\ref{fig:pwi_demo}(e,f), this approximation is shown in Fig.~\ref{fig:pwi_demo}(g). An obvious feature of the exceptional set is the open circular regions, or \textit{cells}, which are periodic non-mixing regions that rotate about an internal elliptic periodic point and are never cut by the cutting lines of the PWI \cite{Park2016a, Park2017, Lynn2019, Smith2017, Scott2001, Scott2003, Goetz2003}. The fractional coverage of the hemispherical shell's area contained in the associated fat fractal is denoted \(\Phi\) and its limiting value is denoted \(\Phi_\infty\) \cite{Park2016a, Park2017, Lynn2019, Smith2017}.  A key point is that while the exceptional set indicates where mixing is possible, mixing is not guaranteed within the exceptional set \cite{Lynn2019}.
 
The purpose of this paper is to investigate invariant ergodic subsets within \(\bar{E}\) that do not intermix. By ergodic set, we mean that every \textit{orbit} (i.e.\ the set of all pre- and post-images of a single point, e.g. \(X = \{M^n(\bm{x}):-\infty<n<\infty\}\) is the orbit of \(\bm{x}\)) from the set either exists in a lower dimensional set (zero measure, e.g., a finite number of periodic points or curves) or fills out the entire ergodic set (full measure); a set that is only partially filled by an orbit is not ergodic, but the subset that is filled may be ergodic on its own \cite{Trinh2014}.  In this case, minimally invariant sets (orbits under the PWI and their closures) with non-zero measure are also ergodic due to the dynamics of PWIs \cite{Umberger1985}. In other words, an ergodic set is one in which, under the PWI map, an orbit passes arbitrarily close to every point within the set infinitely often (Poincar\'e recurrence) and with equal frequency (with respect to an invariant measure) throughout the set. The closure of every orbit that has non-zero measure under the PWI is invariant under the map and defines an ergodic subset of the domain (due to the dynamics of PWIs). If the exceptional set \(\bar{E}\) can be split into multiple, separate ergodic subsets, then there are additional barriers to mixing within the exceptional set. We present new tools to examine invariant subsets in \(\bar{E}\) and definitively determine that \(\bar{E}\) is not always ergodic, and instead contains many ergodic subsets for some PWIs.

Our goal in doing this is to answer two specific questions about the physical nature of solid mixing by cutting-and-shuffling. First, do more cuts necessarily generate more mixing? We demonstrate that this is not the case, and, in fact, variation in the density of cuts at different locations throughout the domain is an indicator of barriers to mixing. Secondly, do non-trivial cutting actions necessarily produce any mixing at all? We demonstrate, by use of a modified recurrence plot, that there exist invariant subsets of cut regions which indicate barriers to mixing, potentially barring global mixing. Both of these results are likely to be of importance for practical applications of solids mixing.

\section{Mathematics of Piecewise Isometries}

Mathematically, a PWI on a metric, measure space \(S\) is a map \(\Phi: \mathcal{P}(S) \to \mathcal{P}(S)\), where \(\mathcal{P}(S)\) denotes the power set of \(S\) (the set of subsets of \(S\)). \(S\) consists of a finite number of \textit{tiles} or \textit{atoms} \(\{P_1,...,P_m\}\) with nonempty interiors that form a tiling of \(S\), and \(\Phi\) consists of isometries on these atoms (metric preserving maps)
\begin{equation*}
\phi_k: P_k \to S, \ 1\le k \le m.
\end{equation*}
Although not a requirement for PWI in general, the PWI studied here is almost everywhere invertible such that \(\{\phi_1(P_1),...,\phi_m(P_m)\}\) is also a tiling of \(S\), i.e.\ \(S = \bigcup_{k=1}^m \phi_k(P_k) = \bigcup_{k=1}^m P_k\), and, as such, both tilings, \(\{P_1,...,P_m\}\) and \(\{\phi_1(P_1),...,\phi_m(P_m)\}\) overlap only on their boundaries, \(\partial P_k\) and \(\partial \phi_k(P_k)\) respectively, which are measure zero. The single PWI map, \(\Phi: \mathcal{P}(S) \to \mathcal{P}(S)\), is the combination of these isometries,
\begin{equation*}
\Phi (\bm{x}) := \phi_k (\bm{x}), \text{ when }  \bm{x} \in P_k, \ 1\le k \le m,
\end{equation*}
which is multivalued and discontinuous on the overlapping boundaries of \(P_k\).
We use \(\mathcal{D}\) to denote the collection of these multi\-valued points (which are discontinuities in the PWI), specifically \(\mathcal{D} = \bigcup_{i \neq j} P_i \cap P_j\). The red cutting line \(\mathcal{D}_1\) (due to the first rotation by \(\alpha\)) and the black cutting line \(\mathcal{D}_2\) (due to the second rotation by \(\beta\)) together form this multi\-valued, zero measure set in Fig.~\ref{fig:pwi_demo}(e,f). In many definitions \cite{Ashwin1997, Goetz2000, Goetz2001, Scott2001, Scott2003, Ashwin2002, Goetz2003, Goetz2004}, PWIs are defined on \(\mathcal{D}\) in such a way as to avoid multivalued mappings, but orbits originating from this set, \(\mathcal{D}\), and nearby are used in this paper to understand the mixing regions adjacent to \(E\). Atoms that share a boundary each retain a copy under the map, which will later be described as each retaining a \textit{side} of the discontinuity.

The multi\-valued map for the hemispherical shell PWI is formally defined as a rotation for each atom independently,
\begin{equation}
M_\ab (\bm{x}) = \{ R_i (\bm{x}) : \bm{x} \in P_i\},
\end{equation}
where the four rotation isometries \(R_i\) for each \(P_i\) are
\begin{eqnarray}
R_1(\bm{x}) &=& R_\beta^{x}R_\alpha^{z}\bm{x},\\
R_2(\bm{x}) &=& R_{\beta }^{x}R_{\alpha+ \pi}^{z}\bm{x},\\
R_3(\bm{x}) &=& R_{\beta + \pi}^{x}R_{\alpha }^{z}\bm{x},\\
R_4(\bm{x}) &=& R_{\beta + \pi}^{x}R_{\alpha + \pi}^{z}\bm{x},
\end{eqnarray}
and \(R_\theta^{a}\) is a rotation about axis \(a\) by angle \(\theta\) \cite{Smith2017,Lynn2019}. Rotations are applied right to left such that rotation about the \(z\)-axis is first and rotation about the \(x\)-axis is second. Additional rotations by \(\pi\) are included where the specified atom crosses the (periodic) equator.

This PWI is orientation preserving (free of reflections) and almost everywhere invertible except for the points that map to the equator after either rotation (atom boundaries) which map to two locations. The pre-image of the points that encounter the equator during the rotation procedure in Fig.~\ref{fig:pwi_demo} are labeled \(\mathcal{D}_1\) (first rotation about the \(z\)-axis) and \(\mathcal{D}_2\) (subsequent rotation about the \(x\)-axis) in Fig.~\ref{fig:pwi_demo}(c). Since \(M \mathcal{D}_2 = \partial S\), the same blue color is assigned to \(\partial S\), \(\mathcal{D}_2\), and their pre- and post-images in later figures to aid visualization. Likewise, the pre- and post-images of \(\mathcal{D}_1\) are red in later figures. Figure~\ref{fig:pwi_demo}(g) shows \(\tilde{E}_{0.001,20\,000}\) (\(N=20\,000\) iterations of cutting lines given half-width \(\varepsilon = 0.001\)), for the \((57^\circ, 57^\circ)\) protocol, an approximation of the exceptional set, \(\bar{E}\), which is the accumulation of all possible pre- and post-images of the cutting lines using this red-blue coloring viewed from the negative \(y\)-axis and flattened to a disk using a Lambert azimuthal equal-area projection \cite{Snyder1987}. This projection preserves the relative areas of features and is used in Sec.~\ref{sec:coloring} to evenly sample across the shell (an evenly spaced grid is projected onto the hemisphere) \cite{Lynn2019}.

The hemispherical shell, \(S\), can be partitioned into distinct invariant sets with respect to the PWI. \(O \subset S\) is the collection of periodic islands, maximally open neighborhoods about central periodic points, which almost always, with the exception of degenerate protocols without periodic islands, has positive Lebesgue measure, i.e.\ 2D area. Periodic islands are never cut by the PWI. In contrast, \(E\) is the countable union of pre- and post-images of the cutting lines, \(\mathcal{D}\), and, as a collection of great circle arcs, has zero Lebesgue measure. \(E\) can contain periodic points, but since \(M\) is treated as a multi\-valued map on \(\mathcal{D}\), it may be that only one of the multiple images of a point in \(E\) is actually periodic. \(\bar{E}\), the closure of \(E\), is called the \textit{exceptional set}. While \(E\) may contain periodic points, the \textit{remainder} of \(\bar{E}\) outside of \(E\), \(E' = \bar{E} \setminus E\), necessarily cannot contain any periodic points \cite{Bruin2003} \footnote{Take periodic orbit \(X\) (such that every \(\bm{x} \in X\) is a periodic point) and the distance metric \(d(\bm{a},B) = \inf_{\bm{b} \in B} d(\bm{a},\bm{b}) \ge 0\) which is the minimum distance between a point \(\bm{a}\) and a set \(B\). If \(d(\bm{x},\mathcal{D}) > 0\) for all \(\bm{x} \in X\), then there exists a \(\delta >0\) such that \(d(\bm{x},\mathcal{D}) > \delta\) for all \(\bm{x} \in X\) and, therefore, there exists a neighborhood around each \(\bm{x} \in X\) constituting a periodic cell with non-zero measure such that \(X \subset O\). If \(d(\bm{x},\mathcal{D}) = 0\) for some \(\bm{x} \in X\), then there is at least one point \(\bm{y} \in X\) such that \(d(\bm{y},\mathcal{D}) = 0\) and \(\bm{y} \in \mathcal{D}\), which implies \(X \subset E\). Therefore, there are no periodic points in \(E'\).}. Since \(E\) has zero measure as a collection of thin arcs, all of the measure associated with a fat fractal \(\bar{E}\) is from \(E'\). These three sets, \(E\), \(E'\) and \(O\) are all invariant under \(M\) (and, of course, \(\bar{E} = E \cup E'\) is trivially invariant as the union of two invariant sets) \footnote{\(E\) is invariant by definition as it is the smallest set containing \(\mathcal{D}\) that is invariant. \(O\) is also invariant as the image of any periodic cell is another cell. \(E' = S \setminus (O \cup E)\) must then also be invariant. If it were not invariant, there would exist \(x \in E'\) such that \(M(x) \notin E'\), i.e.\ \(M(x) \in E \cup O\), which violates the invariance of \(E\) or \(O\).}.
 
In general, \(O\) is not minimally invariant (containing no smaller invariant subset), as each periodic cell can only map to cells of the same size   \cite{Smith2017}. For any given protocol, it is not clear whether \(\bar{E}\) is minimally invariant (containing no smaller invariant subsets within them), but the existence of distinct, positive measure, minimally invariant subsets in \(\bar{E}\) introduces barriers to the mixing induced by the PWI. Orbits in \(\bar{E}\) can be approximated by nearby trajectories in \(E\), but due to roundoff error, no numerically evaluated orbits will ever enter \(E\). 

With this background, we begin by examining an approximation to the natural invariant measure of \(\bar{E}\) and the numerical considerations in representing orbits via color-coding in Sec.~\ref{sec:natural_invariant}. In Sec.~\ref{sec:return_maps} we examine how images of \(\mathcal{D}\) return to \(\mathcal{D}\) using a type of return plot and draw conclusions about invariant subsets within \(\bar{E}\) that can create barriers to mixing. Again, our goal is not to attempt to characterize the infinite range of possible PWIs on the hemispherical shell, which has been done previously \cite{Park2016a, Park2017, Lynn2019, Smith2017}, but instead to develop approaches to understand the mixing characteristics of PWIs in general by considering a few specific hemispherical shell PWIs as examples.

\section{Color-coding invariant subsets} \label{sec:natural_invariant}

\begin{figure}
  \includegraphics[width=0.45\textwidth]{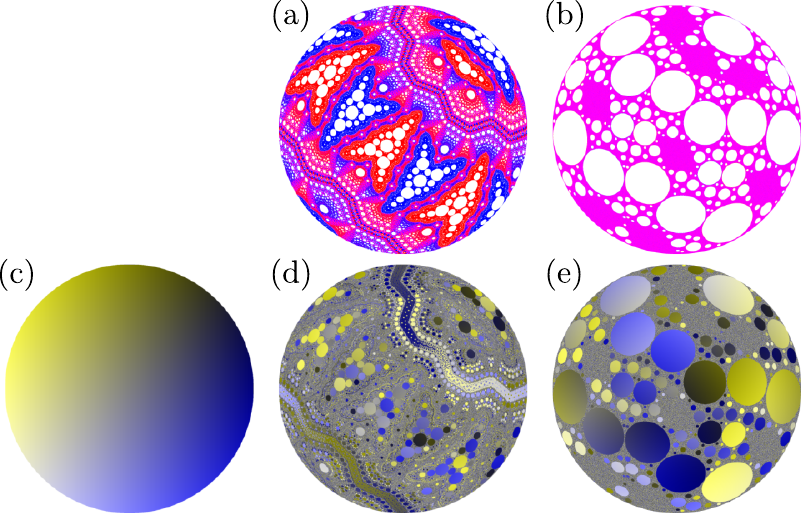}
  \caption{Approximate exceptional sets , \(\tilde{E}_{0.001,2\times10^4}\), for (a) \((45^\circ,45^\circ)\) and (b) \((57^\circ, 32.75^\circ)\) protocols viewed from below using a Lambert equal area projection (see Sec.~\ref{sec:natural_invariant} for explanation of how the hue is computed).
  Corresponding mixing of initial condition in (c) for (d) \((45^\circ,45^\circ)\) and (e) \((57^\circ, 32.75^\circ)\) after \(2\times 10^4\) iterations of \(M_\ab\).
  }
  \label{fig:coverage_compare}
\end{figure}

Previous studies on the hemispherical PWI \cite{Smith2017, Lynn2019} used a correlation \cite{Park2017} between the measure of \(\bar{E}\) (2D area) and the \textit{intensity of segregation} due to Danckwerts \cite{Danckwerts1952}, a measure of mixing, to compare the degree to which two different protocols mix. Some PWIs have ergodic exceptional sets \cite{Kahng2002}, some PWIs produce weak mixing within their exceptional sets \cite{Avila2007}, and other PWIs have exceptional sets containing separate invariant subsets that do not mix between each other \cite{Ashwin2005, Ashwin2018}. In general, it cannot be predicted whether a given PWI will have an ergodic exceptional set, and, as a result, the measure of \(\bar{E}\) is not always an accurate indicator of mixing. Compare, for example, the protocols \((45^\circ, 45^\circ)\) and \((57^\circ, 32.75^\circ)\), shown in Fig.~\ref{fig:coverage_compare}(a,b), respectively, both of which have exceptional sets that cover roughly 41\% of \(S\), i.e.\ \(\Phi_\infty \approx 0.41\) \cite{Park2017, Smith2017, Lynn2019}. For both protocols, the \(N\)-iterate approximation to \(\bar{E}\), \(\tilde{E}_{\varepsilon,N}\), shown in Fig.~\ref{fig:coverage_compare}(a,b) is formed by using \(N = 2\times 10^4\) iterations of fattened cutting lines with width \(2\varepsilon = 0.002\) (\(\varepsilon\) on either side of the line). Using this approach, \(\varepsilon\)-fattened cutting lines subject to the PWI completely cover \(\bar{E}\) in a finite number of iterations \cite{Lynn2019}.

It is illustrative to compare the mixing for these two protocols. The continuously varying initial condition shown in Fig.~\ref{fig:coverage_compare}(c) is partially mixed by both protocols as shown in Figs.~\ref{fig:coverage_compare}(d) and (e). As expected, non-mixing cells do not mix with the rest of the domain, but, surprisingly, the \((45^\circ, 45^\circ)\) protocol in Fig.~\ref{fig:coverage_compare}(d) generates non-circular regions inside of \(\tilde{E}_{\varepsilon,N}\) that do not mix with the rest of the domain. Previous measurements of fractional coverage \cite{Park2017, Lynn2019, Smith2017} count all of \(\bar{E}\) as mixing, yet, there appear to be regions within \(\bar{E}\) that are isolated from one another.  The most obvious example is the zigzag band across the upper right and lower left of Fig.~\ref{fig:coverage_compare}(d). Some parts of this band are collections of periodic cells that travel together, but more interesting is that the parts of this band outside of these cells are completely self contained, indicating a separate invariant subset within \(\bar{E}\). This zigzag band appears to be an interval exchange transform embedded into the PWI \cite{Ashwin2018}. An interval exchange transformation is the one-dimensional version of cutting-and-shuffling in which a line segment (or other one-dimensional object) is split into pieces that are reordered \cite{Yu2016, Novak2009, Katok1980, Hmili2010, Keane1977, Masur1982, Wang2018}.  A closer examination of the \((45^\circ,45^\circ)\) protocol in Sec.~\ref{sec:return_maps} reveals many other self-contained invariant subsets in \(\bar{E}\), each of which represents an isolated mixing region and a barrier to mixing that also indicates \(\bar{E}\) is not ergodic as a whole. In contrast, the \((57^\circ, 32.75^\circ)\) protocol in Fig.~\ref{fig:coverage_compare}(e) produces good mixing (gray regions), except for the unmixed cells corresponding with the white regions in Fig.~\ref{fig:coverage_compare}(b). Thus, this protocol does not appear to have any structure outside of the non-mixing islands, which suggests that there are no smaller invariant subsets in in \(\bar{E}\).

These examples show that although the fractional coverage \(\Phi_\infty\) indicates what fraction of the domain is outside of non-mixing cells, it fails to indicate the degree of mixing within \(\bar{E}\). Hence, the fractional coverage of \(\bar{E}\) only indicates how much mixing is possible, not whether it occurs. In previous work \cite{Park2017, Smith2017, Lynn2019, Scott2003, Ashwin1997a, Goetz2000, Goetz2004}, exceptional sets for the hemispherical PWI, as well as other PWIs, have been constructed without considering the dynamics within them. However, in the process of generating the approximation \(\tilde{E}_{\varepsilon,N}\), the local amount of cutting in the exceptional set, which will be subsequently defined as the \textit{density}, \(\rho\), can be approximated by layering \(\varepsilon\)-fattened (finite width \(2\varepsilon\)) cutting lines on top of one another. Variation in this \(\varepsilon\)-fattened cutting line density reveals the subtle non-mixing structures within \(\bar{E}\), as described in this section.

Two conditions must be true for an exceptional set to be completely ergodic on its own, that is, without separate ergodic subsets. First, since the PWI is composed of isometries (which are non-distorting), the entire ergodic exceptional set must encounter the cutting lines, \(\mathcal{D}\), with equal frequency or density (i.e. the local distribution of cuts is uniform). That is, if one point encounters a cutting line more often than another, the two points must lie on different orbits (and therefore different invariant sets) within the exceptional set. Second, almost any point in an ergodic exceptional set (with the exception of measure zero orbits, e.g. periodic points) necessarily has an orbit dense in the entire exceptional set which includes the initial cutting lines, \(\mathcal{D}\), such that any point in \(\bar{E}\) will have a trajectory that falls arbitrarily close to every point along the entire length of \(\mathcal{D}\), and, in conjunction with the first condition, encounters all of \(\mathcal{D}\) with equal frequency or density. These two conditions motivate the following measurement of the cutting line density.

\subsection{Cutting line density}

Cutting line \textit{density} is an invariant measure of how often a region of the domain is cut. Using the \(\varepsilon\)-fattened cutting line approximation \cite{Lynn2019}, this can be thought of as layering paint along each cutting line such that the paint thickness, when normalized by the number of cuts painted, represents the local density of cutting lines.

Before formally defining the cutting line density \(\rho\), we address some mathematical preliminaries. A mathematical \textit{ball} is used to consider neighborhoods of points (a point plus points in the region surrounding it) such that given a point \(\bm{x}\in S\), the closed ball of radius \(r > 0\), \(B_r(\bm{x}) = S\cap \{\bm{z} \in S : \|\bm{x} - \bm{z}\| \le r\}\), contains all points \(\bm{z}\) within distance \(r\) of point \(\bm{x}\) \footnote{The geodesic distance on the unit sphere is used here. Geodesic distance, \(d_G\), and Euclidean distance, \(d_E\), on the unit sphere are related by \(d_E^2 = \sin^2(d_G) + (1-\cos(d_G))^2\) and are equivalent in the limit as \(d_E\) or \(d_G \to 0\).}. This is generalized to the closed \(r\)-neighborhood of a set \(A\) which is defined as \(B_r(A) = S \cap \bigcup_{\bm{x} \in A} B_r(x) = S \cap \bigcup_{\bm{x} \in A} \{\bm{z} \in S : \|\bm{x} - \bm{z}\| \le r\} \) and likewise is the set of all points \(\bm{z}\) within a distance \(r\) of \(A\). 

The \(\varepsilon\)-fattened cutting lines of the PWI are defined as \(\mathcal{D}_{1,\varepsilon} = B_\varepsilon(\mathcal{D}_1)\) and \(\mathcal{D}_{2,\varepsilon} = B_\varepsilon(\mathcal{D}_2)\) for some small \(\varepsilon > 0\). Note that the two dimensional Lebesgue measure (\(\mathcal{L}^2(\cdot)\), or area) of the overlap is \(\mathcal{L}^2(\mathcal{D}_{1,\varepsilon} \cap \mathcal{D}_{2,\varepsilon}) >0\) for \(\varepsilon > 0\) because there is in fact a small region of overlap at the intersections of the cutting lines \(\mathcal{D}_{1,\varepsilon}\) and \(\mathcal{D}_{2,\varepsilon}\). 
The \(\varepsilon\)-fattened cutting lines are used to define an \(\varepsilon\)-fattened \(E\) \cite{Lynn2019}, which is called \(E_\varepsilon = \bigcup_{n = -\infty}^\infty M^n \mathcal{D}_\varepsilon = B_\varepsilon(E)\) where \(\mathcal{D}_\varepsilon = \mathcal{D}_{1,\varepsilon} \cup  \mathcal{D}_{2,\varepsilon} \), that approximates \(\bar{E}\) (contains \(E\) and its limit points as well as extra points due to the \(\varepsilon\)-fattened cutting lines) since in the limit as \(\varepsilon \to 0\), \(E_\varepsilon \to \bar{E}\) (extra points are reduced and only \(E\) and its limit points remain) \cite{Grebogi1985, Lynn2019}. It is not possible to apply infinite iterations of \(M\) to completely construct \(E_\varepsilon\), so \(N\) iterations are used to find an estimate of \(E_\varepsilon\): \(\tilde{E}_{\varepsilon,N} = \bigcup_{n=0}^N M^n \mathcal{D}_\varepsilon\). There always exists a sufficiently large \(N\) such that the approximation \(\tilde{E}_{\varepsilon,N}\) composed of the \(\varepsilon\)-fattened cutting lines completely covers \(\bar{E}\) \cite{Lynn2019}, and this estimate is increasingly accurate as \(\varepsilon \to 0\) and \(N \to \infty\).

Before defining the cutting line density, we define the fraction, \(F\), of iterates for which \(\bm{x}\) is contained in some set \(A \subset S\) as
\begin{eqnarray}
  F(\bm{x}, A) &=& \lim_{n\to \infty} \frac{\# \{ \bm{x} \in M^i (A) : 0 \le i \le n \}}{n} \label{eq:moving_set} \\
   &=& \lim_{n\to \infty} \frac{\# \{M^{-i} (\bm{x}) \in A : 0 \le i \le n \}}{n}. \label{eq:frac_point}
\end{eqnarray}
where \(\#\) denotes the counting measure such that the numerator is the count of iterates of \(\bm{x}\) in \(A\) \cite{Alligood1996}. The `moving-set' (\(M\) is applied to \(A\)) definition, Eq.~\ref{eq:moving_set}, is equivalent to the more useful (when considering orbits) `moving-point' definition in Eq.~\ref{eq:frac_point}. Equation~\ref{eq:frac_point} does not include limit points arbitrarily close to the measurement set \(A\) which makes it a poor measure for orbits that are arbitrarily close to \(A\). For orbits arbitrarily close, but not contained in \(A\), \(F\) is always 0. To include these limit points, the density of the orbit starting at \(\bm{x}\) in a set \(A\) is defined as
\begin{equation}
	\rho(\bm{x} , A) = \lim_{r \to 0} F(\bm{x}, B_r(A)), \label{eq:density}
\end{equation}
such that limit points are included in this measurement. For attracting maps, this density is identical to the \textit{natural invariant measure} of the set \(A\) under said attracting map \cite{Weisstein, Alligood1996, Ott2002} \footnote{The natural invariant density is a global property (any starting position \(\bm{x}\) will produce the same value) for attracting maps describing the normalized fraction of `time' a typical orbit spends in \(A\) and can be used as a measure on the attractor for attracting maps. Since PWIs are non-attracting maps, this density can be measured throughout the exceptional set and will produce different values for starting points \(\bm{x}\) that are in different invariant sets. The natural invariant density is therefore not a global property of the map (e.g.,  \(\rho(\bm{x},E) = 0\) for all \(\bm{x} \in O\) since \(E\) and \(O\) are distinct invariant sets, \(E \cap O = \emptyset\)).}.

For any \(A \subset S\) it is straightforward to show that \(\rho(\bm{x},A) = \rho(M^n \bm{x},A)\) for any integer \(n\), i.e.\ \(\rho\) is invariant under the map and constant along the orbit of any \(\bm{x}\) in the domain. By extension, if \(C \subset S\) is an ergodic subset, then \(\rho(\bm{z},D)\) will be identical for every \(\bm{z} \in C\) and any subset \(D \subset S\) as a direct result of ergodicity, i.e.\ \(\rho\) is constant throughout an ergodic set.

The sets of interest for measuring the density of the exceptional set are the \(\varepsilon\)-fattened cutting lines that generate the fat fractal. In fact, \(\lim_{\varepsilon \to 0} \rho(\bm{x},\mathcal{D}_\varepsilon)\) is the natural invariant density in \(\mathcal{D}\) \cite{Ott2002, Alligood1996}, but is zero for all non-periodic \(\bm{x}\) since the PWI has no attractors.

Finally, the density of iterates of \(\mathcal{D}_{1,\varepsilon}\) at a point \(\bm{x} \in S\) is \(\rho(\bm{x},\mathcal{D}_{1,\varepsilon})\) (or equivalently, the density of the orbit of \(\bm{x}\) in \(\mathcal{D}_{1,\varepsilon}\)) and the density of \(\mathcal{D}_{2,\varepsilon}\) is likewise \(\rho(\bm{x},\mathcal{D}_{2,\varepsilon})\). As the areas of \(\mathcal{D}_{1,\varepsilon}\) and \(\mathcal{D}_{2,\varepsilon}\) scale linearly with \(\varepsilon\) in the limit of small \(\varepsilon\) and the map is composed of isometries, \(\rho(\bm{x},\mathcal{D}_{1,\varepsilon})\) and \(\rho(\bm{x},\mathcal{D}_{2,\varepsilon})\) are also proportional to \(\varepsilon\) in the limit \(\varepsilon \to 0\). As a consequence, density measurements are normalized by \(\varepsilon\) to more easily compare densities across values of \(\varepsilon\). The error in this measurement is discussed in Appendix~\ref{sec:app_A_errors} in the Supplementary Material. Normalizing in this way gives, for non-periodic orbits, the almost \(\varepsilon\)-invariant quantity \(\rho/\varepsilon\), which, as discussed in Sec.~\ref{sec:density_integrals} and Appendix~\ref{sec:app_B_integral_proof} in the Supplementary Material, is closely connected to the area of the corresponding orbit in the limit as \(\varepsilon \to 0\). 

We consider two representations of this cutting line density. First, a coloring method based on measurements of cutting line density is used to create a colored exceptional set in which the color intensity represents the cutting line density, as described next. Second, a modified recurrence plot, called a \textit{return plot}, that captures localized density throughout the approximate exceptional set is constructed and analyzed in Sec.~\ref{sec:return_maps}, revealing barriers to mixing and confined trajectories within \(\bar{E}\).

\subsection{Color-coding the exceptional set} \label{sec:coloring}

\begin{figure}
  \includegraphics[width=0.45\textwidth]{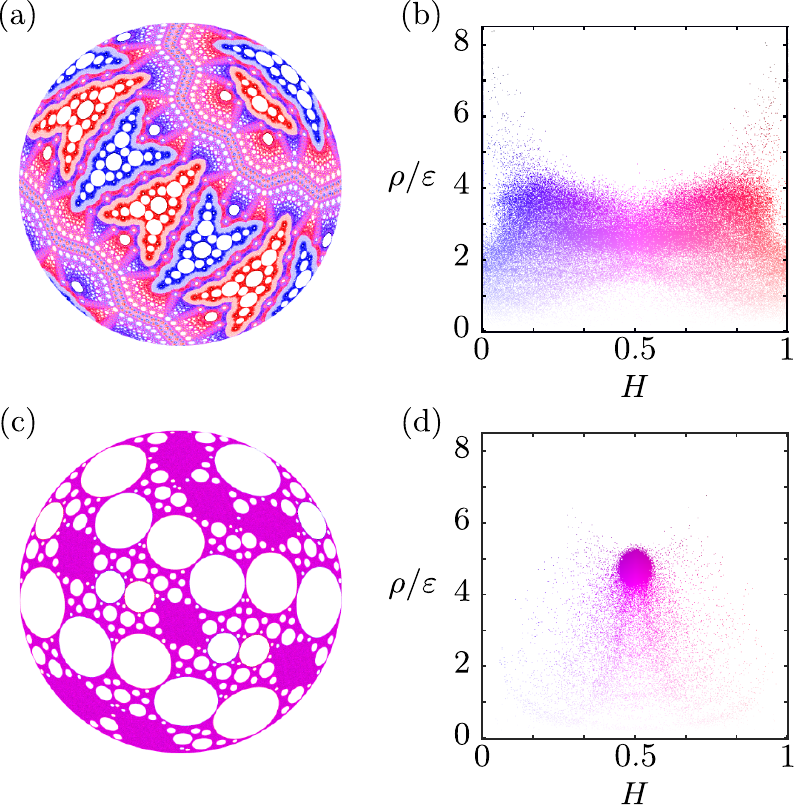}
  \caption{ Area-preserving projection of (a) the \((45^\circ, 45^\circ)\) protocol and (c) the \((57^\circ, 32.75^\circ)\) protocol after \(N = 2\times 10^4\) iterations using a cutting line width of \(\varepsilon = 0.001\) as viewed from below. (a,c) Both hue (relative density) and lightness (density) of cutting lines. (b,d) Scatter plot of hue, \(H\), and cutting line density, \(\rho/\varepsilon\), for colors shown in (a,c) respectively.}
  \label{fig:45combined}
\end{figure}

With this background, we consider the cutting line density for an example protocol \((45^\circ, 45^\circ)\) after \(2 \times 10^4\) iterations. Applying the PWI repeatedly, the cutting lines, \(\mathcal{D}_1\) and \(\mathcal{D}_2\), are cut-and-shuffled around the domain forming a part of \(E\), which has cutting lines overlapping only at individual points. With \(\varepsilon\)-fattened cutting lines, cutting line intersections are small patches instead of individual points, that is \(\mathcal{L}^2(\mathcal{D}_\varepsilon \cap M^i \mathcal{D}_\varepsilon) > 0\). This is what contributes to density (individual points do not contribute to density). These small overlapping areas are categorized in two ways to make the overall coloring shown in Fig.~\ref{fig:45combined}(a): the total density (of both fattened cutting lines) is used for lightness and the relative density of the two different cutting lines is used for hue. The hue measurement alone is shown in Fig.~\ref{fig:coverage_compare}(a). The different hues in Fig.~\ref{fig:coverage_compare}(a) demonstrate that some parts of the domain are cut mostly (or even exclusively) by \(\mathcal{D}_1\) (red), others are cut mostly by \(\mathcal{D}_2\) (blue), and still others are cut by both \(\mathcal{D}_1\) and \(\mathcal{D}_2\) (magenta).

When the absolute density is also included as a lightness with the hue, the image in Fig.~\ref{fig:45combined}(a) is produced. White areas show cells in which there are no cutting lines at all (and thus no cutting line density). The lightness of color in the \(\varepsilon\)-fattened exceptional set in Fig.~\ref{fig:45combined}(a) reflects the total cutting line density, \(\rho(\bm{x},\mathcal{D}_\varepsilon)/\varepsilon\). Vivid color represents a high cutting line density (many cuts), whereas pale color indicates low cutting line density (few cuts).

Thus, in visualizing the exceptional set, two different color dimensions are used: lightness based on the total cutting line density, and hue based on the relative densities of the original, red-blue cutting lines. An HSL [hue, saturation (constant), lightness] color space is used to show both dimensions together. The lightness at a point \(\bm{x} \in S\) is defined as the total density of the two cutting lines,
\begin{equation}
	L(\bm{x}) = \rho(\bm{x},\mathcal{D}_{1,\varepsilon}) + \rho(\bm{x},\mathcal{D}_{2,\varepsilon}). \label{eq:lightness}
\end{equation}
As noted earlier, since \( \mathcal{L}^2 (\mathcal{D}_{1,\varepsilon}\cap \mathcal{D}_{2,\varepsilon}) > 0\) for \(\varepsilon > 0\), this lightness is not exactly density of \(\mathcal{D}_\varepsilon\), \(L(\bm{x}) \neq \rho(\bm{x},\mathcal{D}_\varepsilon)\), because \(L(\bm{x})\) double-counts the region at the intersection of \(\varepsilon\)-fattened cutting lines. However, since \( \mathcal{L}^2 (\mathcal{D}_{1,\varepsilon}\cap \mathcal{D}_{2,\varepsilon})\) is roughly equal to \(2\varepsilon^2[\csc(\alpha)\sec(\alpha) + \csc(\beta)\sec(\beta)]\) for small \(\varepsilon\), the two measurements are similar, \(L(\bm{x}) \sim \rho(\bm{x},\mathcal{D}_\varepsilon)\), in the limit of small \(\varepsilon\) for \(\alpha, \beta \gg 0\) (see Appendix~\ref{sec:app_A_errors}  in the Supplementary Material for more discussion on the error incurred by this overlap). Hence, we use the symbols \(\rho\) and \(L\) interchangeably to mean \(L(\bm{x}) \approx \rho(\bm{x},\mathcal{D}_\varepsilon)\). For small values of \(\alpha\) or \(\beta\), this similarity breaks down, because the overlap of \(\mathcal{D}_1\) and \(\mathcal{D}_2\) is large. For an ergodic set, the lightness will be uniform throughout since every point returns to the cutting lines with equal frequency.

The second color dimension is the hue. Hue is measured as a fraction of interactions with one of the two cutting lines out of the total interactions with both cutting lines,
\begin{equation}
	H(\bm{x}) = \frac{\rho(\bm{x},\mathcal{D}_{1,\varepsilon})}{\rho(\bm{x},\mathcal{D}_{1,\varepsilon}) + \rho(\bm{x},\mathcal{D}_{2,\varepsilon})}, \label{eq:hue}
\end{equation}
which is undefined when both \(\rho(\bm{x},\mathcal{D}_{1,\varepsilon})\) and \(\rho(\bm{x},\mathcal{D}_{2,\varepsilon})\) are zero (i.e.\  inside periodic cells). The range of \(H(\bm{x})\) is \([0,1]\), from blue to red. For an ergodic set, the hue is uniform throughout the set since every point encounters either cutting line with equal frequency.

As mentioned above, an ergodic subset \(B \subset S\) has the same density, \(\rho(\bm{x},C)\)  for all \(\bm{x}\in B\) for any sample set \(C\). Equal densities in an ergodic set imply the same coloring according to this visualization method. Although identical color cannot be used to verify that a subset is indeed ergodic, differences in color can be used to show that subsets are necessarily not ergodic. The scatter plot of hue and lightness in Fig.~\ref{fig:45combined}(b) clearly shows that the \((45^\circ,45^\circ)\) protocol is not ergodic since color is not uniform in either hue or lightness. Similar color-coding has been used by Ashwin and Goetz \cite{Ashwin2005} to color-code periodic cells based on their return frequency to a particular atom of the PWI, \(\rho(\bm{x},P_i)\).

The hue of the \((45^\circ,45^\circ)\) protocol, calculated using Eq.~\ref{eq:hue} and shown in Fig.~\ref{fig:coverage_compare}(a), reveals large arrowhead like structures that are exclusively red or blue that do not mix with the rest of the domain. From the lightness, combined with the hue in Fig.~\ref{fig:45combined}(a), it is clear that there are even regions within these arrowheads that are not mixing within themselves, notably their lighter borders have a low cutting line density and do not mix with their darker cores, which have a different (higher) cutting line density.
Furthermore, recall that Fig.~\ref{fig:coverage_compare}(d) shows that the \((45^\circ, 45^\circ)\) protocol has invariant subsets within \(\tilde{E}\) that hinder mixing, such as the zigzag bands across the upper right and lower left. When using the normalized density in Fig.~\ref{fig:45combined}(a), these zigzag bands are a light magenta indicating low cutting line density bands that separate darker regions of high cutting line density in the exceptional set from one another.

For a protocol to have an ergodic \(\bar{E}\), there are conditions on what the coloring of \(\tilde{E}_{\varepsilon,N}\), the approximate exceptional set, can be. Ergodicity implies that dynamics within \(\bar{E}\) should have no bias in proximity to the two cutting lines. Since the average color of \(\mathcal{D}\) is a magenta halfway between red and blue, the resulting stacked cutting lines should be this same average hue of magenta (\(H(\bm{x}) = 0.5, \ \bm{x} \in \bar{E}\)). Further, since \(M\) is a PWI and therefore area preserving, all regions within \(\bar{E}\) should pass within \(\varepsilon\) of cutting lines with an equal frequency. If \(\bar{E}\) is ergodic, \(\tilde{E}_{\varepsilon,N}\) should be a solid magenta color without variation in intensity everywhere inside of \(\bar{E}\) (periodic cells will have single colored boundaries resulting from an intersection between \(E_\varepsilon\) and \(O\)). This does not mean that a single magenta-colored \(\tilde{E}_{\varepsilon,N}\) implies ergodicity, but it does mean that if any portion of \(\tilde{E}_{\varepsilon,N} \cap \bar{E}\) is not a single magenta color, the protocol definitively produces a non-ergodic exceptional set. The distribution of colors using hue and lightness for the \((45^\circ,45^\circ)\) protocol in Fig.~\ref{fig:45combined}(b) indicates that neither hue nor cutting line density are uniform throughout the exceptional set, each varying substantially, and therefore the exceptional set for the \((45^\circ,45^\circ)\) protocol is not ergodic.

The differences in \(\tilde{E}_{\varepsilon,N}\) for the \((45^\circ, 45^\circ)\) protocol shown in Fig.~\ref{fig:45combined}(a) and the \((57^\circ, 32.75^\circ )\) protocol shown in Fig.~\ref{fig:45combined}(c) are obvious. The lightness of cutting line intersections in Fig.~\ref{fig:45combined}(c) shows no distinguishable pattern and is nearly uniform in density. This does not imply that \(\bar{E}\) is ergodic, but it does not rule out ergodicity as it did for \((45^\circ, 45^\circ)\). The hue of overlapped cutting lines in Fig.~\ref{fig:coverage_compare}(b) and \ref{fig:45combined}(c) also shows no discernible pattern within \(\bar{E}\), showing the entire hemispherical shell is uniformly cut by both \(\mathcal{D}_1\) and \(\mathcal{D}_2\) throughout except, of course, in the cells where there are no cuts. Although it cannot necessarily be concluded that \(\bar{E}\) is ergodic, it is likely that \(\bar{E}\) for the \((57^\circ, 32.75^\circ )\) protocol has fewer ergodic subsets than for \((45^\circ, 45^\circ)\). We conclude that mixing within the exceptional set in the \((57^\circ, 32.75^\circ )\) protocol is quite different from that for the \((45^\circ, 45^\circ)\) protocol, despite both exceptional sets occupying nearly the same fraction of the hemisphere, due to the apparent absence of separate invariant subsets and the barriers to mixing they create within the exceptional set in the former. This is emphasized by the tight distribution of hue and cutting line density in the scatter plot in Fig.~\ref{fig:45combined}(d).

\begin{figure}
  \includegraphics[width=0.45\textwidth]{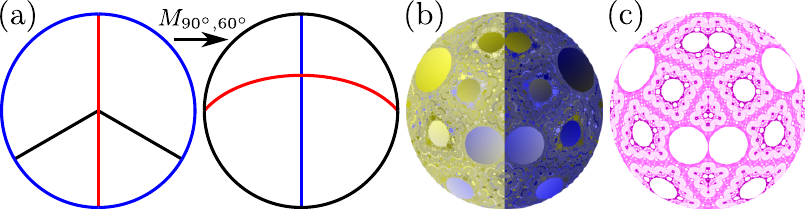}
  \caption{(a) The \((90^\circ, 60^\circ)\) protocol shown from below using an equal-area projection to a disk. (b) Using the initial condition in Fig.~\ref{fig:coverage_compare}(c), mixed state after \(2\times 10^4\) iterations. (c) Hue and lightness visualization of cutting lines for \(\varepsilon =0.001\) after \(2\times 10^4\) iterations of the PWI.
  }
  \label{fig:906090_PWI}
\end{figure}

Unfortunately, neither hue, \(H(\cdot)\), nor density, \(\rho/\varepsilon\), as defined here completely identifies invariant subsets in \(\bar{E}\). To illustrate this, consider the \((90^\circ, 60^\circ)\) protocol shown in Fig.~\ref{fig:906090_PWI}. This protocol has coincident red and blue cutting lines, shown in Fig.~\ref{fig:906090_PWI}(a) with a red cutting line down the vertical diameter for the initial cutting lines and a blue cutting line down the vertical diameter. This creates poor mixing as shown in Fig.~\ref{fig:906090_PWI}(b), but this also makes separation of cutting lines by hue as defined here impossible as shown in Fig.~\ref{fig:906090_PWI}(c), since every red cutting line maps to a blue cutting line and \(\rho(\bm{x},\mathcal{D}_{1,\varepsilon}) = \rho(\bm{x},\mathcal{D}_{2,\varepsilon})\). Notice that all points in \(\tilde{E}_{\varepsilon,N}\) in Fig.~\ref{fig:906090_PWI}(c) have the same magenta hue, and the only variation is due to changes in cutting line density. This protocol produces a left-right barrier to mixing due to the coincident cutting lines \cite{Lynn2019} that is obvious in Fig.~\ref{fig:906090_PWI}(b) but not reflected in the hue-lightness visualization, Fig.~\ref{fig:906090_PWI}(c). In order to differentiate invariant subsets in protocols like this, more information about orbits within \(\bar{E}\) is needed, as described in Sec.~\ref{sec:return_maps}.

\section{Categorizing invariant subsets using recurrence} \label{sec:return_maps}

The coloring scheme in Sec.~\ref{sec:coloring} based on different colors for cutting lines \(\mathcal{D}_1\) and \(\mathcal{D}_2\) is useful for understanding which cutting lines are responsible for different parts of \(\bar{E}\), but this red-blue partitioning is somewhat coarse. Instead of considering just \(\mathcal{D}_1\) and \(\mathcal{D}_2\) separately (by color), the cutting lines can be decomposed into smaller pieces and ultimately into a continuum of points with densities that indicate their interactions.

Every intersection point of iterated cutting lines, \(M^n \mathcal{D}\), with the original cutting lines, \(\mathcal{D}\), is a recurrent point from \(\mathcal{D}\) to \(\mathcal{D}\). Points \(\bm{x} \in \mathcal{D}\) and \(M^n(\bm{x})\) are called recurrent if \(M^n(\bm{x}) \in \mathcal{D}\). Each recurrent point necessarily occupies the same orbit under the PWI. For most protocols, these recurrent points appear to be dense in \(\mathcal{D}\), but there are only a countably infinite number of them \footnote{Recurrent points can only arise due to intersected cutting lines of which there are only ever a countably infinite number. That is, although there are an infinite number of cutting lines in \(E\), they can be ordered according to the iteration of the map that generates them.}, meaning almost all of \(\mathcal{D}\) is not exactly recurrent. This concept of recurrence can be expanded to the \(\varepsilon\)-fattened cutting lines to approximate the exact return structure of the cutting lines. Orbits in \(\bar{E}\) can then be grouped by which points in \(\mathcal{D}\) produce the same orbit. This allows invariant subsets in \(\bar{E}\) to be uniquely identified by the segments of \(\mathcal{D}\) they intersect.

\subsection{Parameterization of cutting lines} \label{sec:paramet}

\begin{figure}
  \includegraphics[width=0.45\textwidth]{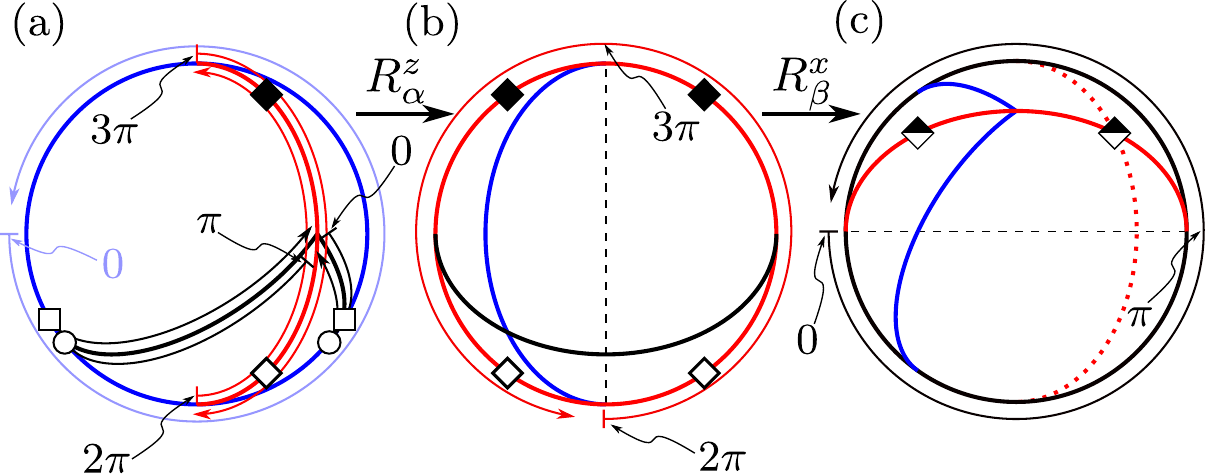}
  \caption{The \((45^\circ,45^\circ)\) PWI's cutting lines' (a) initial position, (b) position after rotation about the \(z\)-axis by \(\alpha\), and (c) position after both rotations; shown orthographically from the negative \(y\)-axis. The black cutting line (and blue boundary) in (a) is parameterized by \(\theta \in [0,2\pi)\) in (c). The red cutting line is parameterized by \(\theta \in [2\pi,4\pi)\) in (b). The parameterization of the initial black cutting line in (a) is broken by the equator and the two pairs of matching points for the parameterization are labeled with circles and squares respectively. The upper right half-filled diamond in (c) is a recurrent point lying on the initial red cutting line (dotted red) and is the combined post-image of the filled and unfilled diamonds in (a).
  }
  \label{fig:parameterization}
\end{figure}

In order to separate orbits in \(\bar{E}\) by their intersections with \(\mathcal{D}\), every point in \(\mathcal{D}\) needs a unique identifier. For the hemispherical shell PWI, atom boundaries that make up \(\mathcal{D}\) are great circle arcs with finite length. This allows a simple parameterization along atom boundaries to uniquely identify each element of \(\mathcal{D}\). For almost all protocols, it is sufficient to parameterize along the length of just the red and black lines in Fig.~\ref{fig:pwi_demo}(e) or just the red and blue lines in Fig.~\ref{fig:pwi_demo}(f) whose total length is \(2\pi\). However, some protocols, such as \((90^\circ, 60^\circ)\) shown in Fig.~\ref{fig:906090_PWI}(a), have a pronounced ``sidedness'' associated with the coincident cutting lines. This is evident in the separation between the two sides of the vertical mid-line, where the cutting lines are coincident, in Fig.~\ref{fig:906090_PWI}(c), indicating that each side of the same cutting line interacts only with its own side of the domain. Because of this, a parameterization should account for the different atoms on either side of the cutting lines. To be clear, the two `sides' of a cutting line are simply the two sets of dynamics assigned to the cutting line, one for each neighboring atom \(P_k\).  This sidedness is also prominent for any protocol that produces polygonal cells whose sides return exactly to \(\mathcal{D}\), including those that produce polygonal tilings, which are composed entirely of periodic segments \cite{Smith2017, Lynn2019}.

The \((45^\circ,45^\circ)\) protocol is used as an example to demonstrate the parameterization of cutting lines, as shown in Fig.~\ref{fig:parameterization}(a). To parameterize along the initial red and black cutting lines in a manner that accounts for the sidedness of the cutting lines, the parameterization is specified based on the orientation, with respect to the rotation procedure in Fig.~\ref{fig:pwi_demo}, that places the cutting line on the hemispherical shell equator. Cutting lines are not multi\-valued when they lie on the equator, which makes measurements accounting for sidedness natural.  Working backwards from the final orientation of the cutting lines after one iteration in Fig.~\ref{fig:parameterization}(c), the black cutting line is parameterized first from \(\theta \in [0,2\pi)\) starting at a point coinciding with the second rotation axis [outer black circle in Fig.~\ref{fig:parameterization}(c)] such that inverting the rotation about the second axis collapses this parameterization onto the black cutting line in Fig.~\ref{fig:parameterization}(a,b). This collapse is why, for protocols without explicitly coincident cutting lines, a sided parameterization often leads to redundant parameter values. The red cutting line is treated similarly and is parameterized according to its orientation when it lies on the equator [outer red circle in Fig.~\ref{fig:parameterization}(b)] from \(\theta \in [2\pi, 4\pi)\) which similarly collapses to the single cutting line in Fig.~\ref{fig:parameterization}(a,c). This parameterization defines the arc-length parameterization of \(\mathcal{D}\) which is referred to as \(\bm{s}(\theta): [0,4\pi) \to \mathcal{D}\). Note that the equatorial circular arc \(\theta \in [0, 2\pi)\) in Fig.~\ref{fig:parameterization}(c) collapses to the two sides of the black cutting line in Fig.~\ref{fig:parameterization}(a) with two matching points in the parameterization (broken by the PWI) labeled by a pair of circles or squares, and that the equatorial circular arc \(\theta \in [2 \pi, 4\pi)\) in Fig.~\ref{fig:parameterization}(b) collapses to the two sides of the red cutting line in Fig.~\ref{fig:parameterization}(a). This parameterization is not unique as any parameterization could be chosen for \(\mathcal{D}\), but this one is selected for its connection to the single-valued orientation of the original cutting lines and because it is an arc-length parameterization.

The necessity for a sided or single-valued parameterization is made clear by examining the points along \(\mathcal{D}\) labeled by filled and unfilled diamonds in Fig.~\ref{fig:parameterization}. These two points in Fig.~\ref{fig:parameterization}(a) map to two opposite points on the equator in Fig.~\ref{fig:parameterization}(b) resulting in four total points. The sided parameterization allows a distinction between these split points to be made before the split occurs. Using a single valued, or sided, parameterization allows the tracking of all four points even when two of them have collapsed onto one another as in Fig.~\ref{fig:parameterization}(a,c). Note that the upper right point in Fig.~\ref{fig:parameterization}(c) has returned to the original cutting line (the dotted red line) such that the sided parameterization breaks down for the next iteration at this point (since it is multivalued for yet another iteration). Despite this breakdown at a single point, the sided parameterization provides additional information about how the cutting lines move subject to the PWI mapping.

\subsection{Return plots}

The cutting line parameterization provides a mechanism for grouping recurrent points in \(\mathcal{D}\), that is, points that start in \(\mathcal{D}\) and return back to \(\mathcal{D}\). Orbits may be trapped within an invariant set, and by separating out points that are confined to certain invariant sets, barriers to mixing within cut regions can be located. In the physical tumbler, orbits confined to cells do not mix with the rest of the domain \cite{Zaman2018}. Other barriers to mixing also seem to exist that, despite being cut by the map, produce confined orbits that do not actually contribute to the overall mixing within the exceptional set.

For example, one recurrent point for \((45^\circ,45^\circ)\) after one application of the PWI is indicated with a half-filled diamond in the upper right of Fig.~\ref{fig:parameterization}(c), but there are recurrent points wherever cutting lines intersect, including the boundary (which is always treated using the same parameterization as the initial black cutting line). Since recurrent points are marked by cutting line intersections, an ``intersection plot'' or ``return plot,'' similar to a recurrence plot for time series \cite{Marwan2009, Eckmann1987, Marwan2007}, can be made where the horizontal and vertical axes are both \(\theta \in [0,4\pi)\). If the point specified by \(\bm{s}(\theta_1)\) returns to \(\mathcal{D}\) at a point specified by \(\bm{s}(\theta_2)\), then a return plot will be populated at both \((\theta_1,\theta_2)\) and \((\theta_2,\theta_1)\), since recurrence is symmetric. This also means that the entire \(\theta_1 = \theta_2\) line will be populated, since every point in \(\mathcal{D}\) is trivially recurrent upon itself. Exact intersections (up to machine precision) can be computed to create an exact intersection map, which is essentially a first return plot for \(\mathcal{D}\), by applying \(M\) to segments of \(\mathcal{D}\) and finding the exact locations of their intersection with other parts of \(\mathcal{D}\), which breaks each segment into smaller segments, using the parameterization in Fig.~\ref{fig:parameterization} to exactly identify points.

Exact recurrence measured in this way has shortcomings. First, the plot will only ever be sparse since cutting line intersections are countable. Additionally, due to the multivalued nature of the discontinuity, any sort of sidedness of cutting lines (neighboring atoms having discontinuous dynamics at the cutting line) cannot be detected computationally using exact recurrence. This is because a segment returning exactly to the original cutting line, either in whole or at its endpoint will not be cut (an example where segments return exactly at their endpoints is the \((90^\circ,60^\circ)\) protocol shown in Fig.~\ref{fig:906090_PWI}) and numerical comparison of two floating point zeros (used to detect intersections) is impossible. Additionally, as cutting line segments are continually split into smaller and smaller segments, the number of cutting line segments to be tracked (in order to find intersections) appears to grow almost linearly with the number of iterations, \(N\), for protocols that do not result in polygonal tilings \footnote{This is a preliminary result that is backed up by evidence from the described method not detailed here, but beyond the scope of this work.}. This results in a constantly growing memory requirement that quickly makes the process of making an exact return plot intractable.

\begin{figure}
  \includegraphics[width=0.45\textwidth]{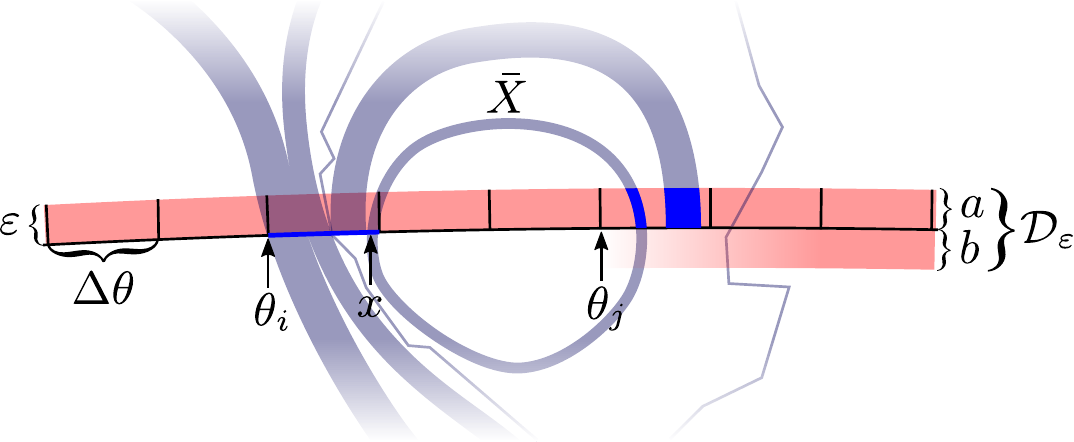}
	\caption{Schematic diagram of a fat, two sided (\(a\) and \(b\)) cutting line in red with example orbits in light gray. The cutting line is split into bins of length \(\Delta\theta\) in which local orbit density is measured. These example orbits are not representative of all possible orbits in an actual PWI.
	}
	\label{fig:cutting_line_assembly}
\end{figure}

Instead of an exact return plot, an approximate return plot is used based on \(\varepsilon\)-fattened cutting lines. Points along the cutting lines are iterated using the PWI, and, if they return to \(\mathcal{D}_\varepsilon\), they approximately return to the closest point in \(\mathcal{D}\). The cutting lines are split into a finite number, \(T\), of bins \([\theta_i,\theta_{i} + \Delta\theta)\) where \(\theta_i, \theta_i + \Delta \theta \in [0,4\pi]\), and each return to the \(\varepsilon\)-fattened cutting lines is placed in a bin within a 2D grid of starting bins and return bins. Figure~\ref{fig:cutting_line_assembly} illustrates a sided, \(\varepsilon\)-fattened cutting line (horizontal with red shading to show the \(\varepsilon\)-fattening) to demonstrate how the cutting line intersects with some example orbits. The cutting line is segmented into \(\Delta \theta\)-long bins in which the local density of an orbit can be measured.

To create the return plot, points are initially placed along \(\mathcal{D}\) in the \(i\)-bins \(\bm{s}([\theta_i,\theta_i + \Delta\theta))\), which are assigned to the horizontal axis of the return plot. These points then increment the corresponding return \(j\)-bin, \([\theta_{j}, \theta_{j} + \Delta\theta)\), along the vertical axis when they return within \(\varepsilon\) of the cutting line as they are iterated by the PWI. This localized return density, called \(\rho_{ij}\), represents how often a point starting in the bin \([\theta_i,\theta_{i} + \Delta\theta)\) returns to the bin \([\theta_j, \theta_{j} + \Delta\theta)\). For example, point \(x\) in Fig.~\ref{fig:cutting_line_assembly} will fill out the orbit \(\bar{X}\) under the PWI, and every return to the \(\theta_j\) bin (dark blue region of \(\bar{X}\) on the right) updates the \(\rho_{ij}\) measurement and indicates that the two bins are connected by an orbit. This can also be thought of as the fractional area of orbits from the \(\theta_i\) bin that intersects the \(\theta_j\) bin as shown in Fig.~\ref{fig:cutting_line_assembly}. Several points to the left of point \(\bm{x}\) in bin \(\theta_i\) also intersect bin \(\theta_j\) following a different orbit, while other points further to the left of point \(\bm{x}\) in bin \(\theta_i\) follow orbits that do not intersect bin \(\theta_j\). Also recall that the \(\varepsilon\)-fattened cutting line has two sides; points on opposite sides of the cutting line may follow different orbits and intersect different bins. Because of this, initial points are placed a small distance \(\delta \ll \varepsilon\) from the cutting line to simplify this distinction in the computational analysis.

Because this local density links a cutting line segment and an \(\varepsilon\)-wide bin, there is an inherent, but subtle, asymmetry in the resulting local density (only disappearing in the \(\varepsilon \to 0\) limit). The local density is written mathematically as the cumulative fraction of all orbits from the \(\theta_i\) bin that intersect the \(\theta_j\) bin,
\begin{equation}
\rho_{ij} = \int_{\theta_i}^{\theta_i + \Delta\theta}\rho(\bm{s}(t),\{\varepsilon  \times [\theta_j,\theta_j + \Delta\theta)\})   \ dt,
\end{equation}
where \(\{\varepsilon  \times [\theta_j,\theta_j + \Delta\theta)\}\) is the \(\theta_j\) bin constructed by the \(\varepsilon\) expansion (to the left in the direction of increasing \(\theta\)) of the cutting line segment making up the entire \(\theta_j\) bin as shown in Fig.~\ref{fig:cutting_line_assembly}.

\begin{figure*}
  \includegraphics[width=0.95\textwidth]{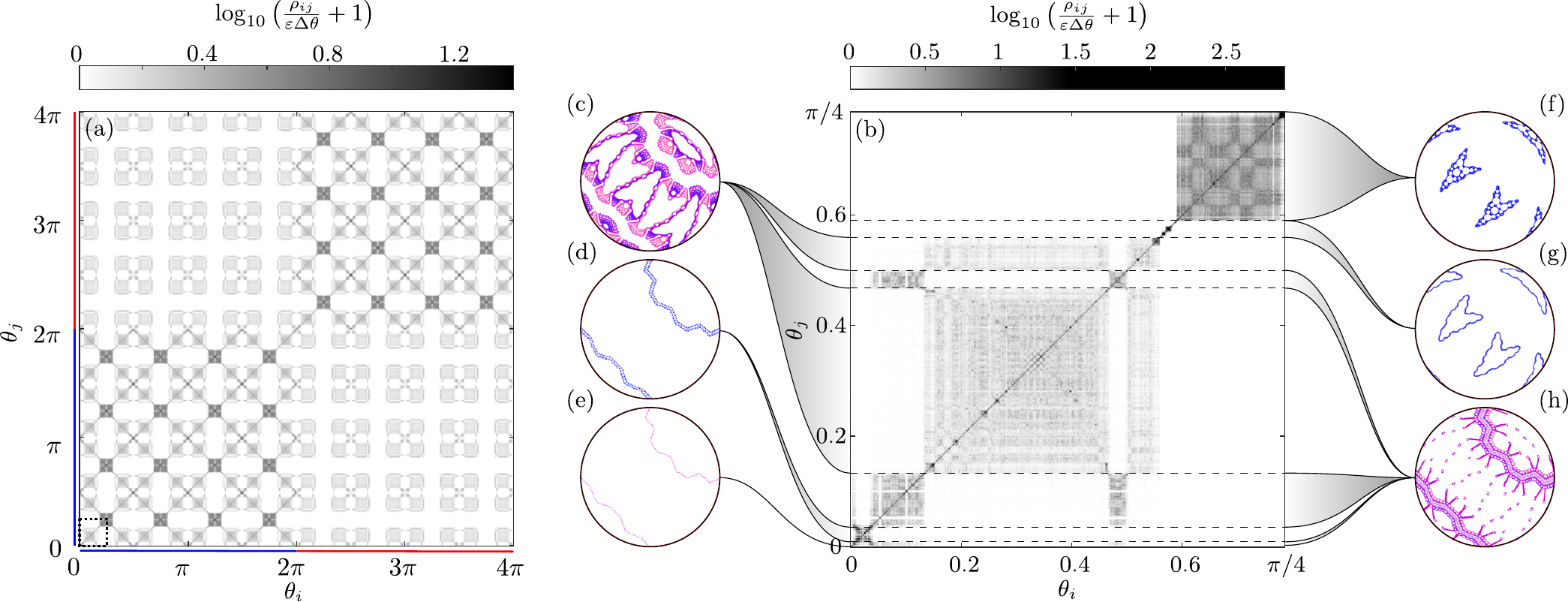}
  \caption{
    (a) Approximate return plot for the \((45^\circ, 45^\circ)\) protocol as a 2D histogram of \(\rho_{ij}\), return to bin \(i\) from bin \(j\); \(\Delta\theta = \pi\times 10^{-3}\). (b) Detail of (a) from the dashed black box in the lower left; \(\Delta\theta = \pi/16\times 10^{-3}\). (a,b) \(2 \times 10^6\) points are initialized \(\delta = 10^{-9}\) from the cutting line, iterated \(N = 2\times 10^6\) times using a cutting line width \(\varepsilon = 10^{-5}\). The \(\theta \in [0,2\pi)\) range represents \(\mathcal{D}_2\) (blue/black cutting line) and the \(\theta \in [2\pi,4\pi)\) range represents \(\mathcal{D}_1\) (red cutting line). (c-h) Various regions of the exceptional set corresponding to portions of the return plot.
  }
  \label{fig:45_return_and_detail}
\end{figure*}

Figure~\ref{fig:45_return_and_detail}(a) shows the resulting return plot for the \((45^\circ,45^\circ)\) protocol. We plot the value of \(\rho_{ij}\) normalized by both the cutting line half-width, \(\varepsilon\), and the bin length, \(\Delta\theta\). Orbits where points return frequently to a small number of bins (possibly localized near the cutting line) result in highly skewed local density \(\rho_{ij}\) values [e.g.,  black dots at the centers of the dark squares in Fig.~\ref{fig:45_return_and_detail}(a)]. The skewed data is transformed to \(\log_{10}\left(\frac{\rho_{ij}}{\varepsilon \Delta \theta} + 1 \right)\) to better reveal the less frequent return regions, which are lighter (due to orbits occupying more space away from the cutting lines). Note that the earlier measurement of hue used to color the exceptional set in Sec.~\ref{sec:natural_invariant} is exactly this measurement with only two bins: one bin for all the points along the red cutting line, \([0,2\pi)\), and one bin for all the points along the blue cutting line, \([2\pi,4\pi)\), indicated on the axes of Fig.~\ref{fig:45_return_and_detail}(a) with adjacent blue and red lines.

The return plot is reminiscent of the adjacency matrix of a weighted graph, with each entry relating the weight of connectivity between two infinitesimal nodes along the cutting lines \footnote{Any introductory text on graph theory (e.g.,  \cite{Ray2013}) should include a discussion of weighted graphs and shortest path problems as well as an algorithm due to Dijkstra for computing the shortest path \cite{Dijkstra1959}.}. The weights in a traditional weighted graph adjacency matrix indicate the difficulty of traveling between nodes. Here, the return plot indicates the reciprocal of this difficulty, i.e.\ the ease of travel between two `nodes,' a value of zero indicating travel is impossible. Questions about the shortest orbit given some diffusion on the scale of \(\varepsilon\), which is built into the return plot, between any two parts of the cutting line can be answered. A measurement of a time-to-\(\varepsilon\)-connectivity between regions of the cutting line, which comes from this perspective, could indicate the rate of mixing within the system, but these analyses are beyond the scope of this study.

For Fig.~\ref{fig:45_return_and_detail}(a), each \([\theta_i,\theta_{i} + \Delta\theta )\) bin is seeded with \(10^3\) points uniformly arranged along the bin. In order to maintain sidedness, the initial points are placed a distance \(\delta \ll \varepsilon\) away from either side of the cutting line. The error in position after \(N\) iterations for the double precision numbers grows like \(N \times (2 \times 10^{-16})\) where \(2 \times 10^{-16}\) is roughly machine precision (sometimes called machine epsilon). As a direct result of this error, it is not possible to place points exactly on a cutting line. The plots shown in the following figures use \(N = 2\times 10^6\) iterations such that the error in position after all iterations is roughly \(4 \times 10^{-10}\). Setting \(\delta\) smaller than this value allows for points to potentially switch sides of the initial cutting line at some point during the application of the PWI, so \(\delta = 10^{-9}\) is used for the following figures. See Appendix~\ref{sec:app_A_errors} in the Supplementary Material for further discussion of the error introduced and mitigated by \(\delta\). Since sidedness often adds a symmetry to protocols without coincident cutting lines such as \((45^\circ,45^\circ)\), the approximate return plot in Fig.~\ref{fig:45_return_and_detail}(a) is symmetrical across \(\theta = \pi\) within the \([0, 2\pi)\) range and across \(\theta = 3\pi\) within the \([2\pi, 4\pi)\) range in both axes. There is also a symmetry across \(\theta_i = \theta_j\), because an orbit connecting \(\theta_i\) to \(\theta_j\) also connects \(\theta_j\) to \(\theta_i\) (there is a subtle asymmetry when using \(\varepsilon\)-fat cutting lines that is invisible in Fig.~\ref{fig:45_return_and_detail}(a); see Appendix~\ref{sec:app_A_errors} in the Supplementary Material).

How the orbits throughout the domain are connected with one another directly affects mixing. Figure~\ref{fig:45_return_and_detail}(a) shows approximately where \(\bar{E}\) is connected and, likewise, where it is disconnected. If \(\rho_{ij}\) is non-zero (gray or black), there exists at least one orbit in \(\bar{E}\) that begins in \([\theta_i, \theta_i + \Delta\theta)\) along the cutting line and either intersects or passes within \(\varepsilon\) of the \([\theta_j , \theta_j + \Delta\theta)\) bin. If \(\rho_{ij}\) is zero, white in Fig.~\ref{fig:45_return_and_detail}(a), then no such orbit exists and bins \(i\) and \(j\) are disconnected from one another such that the orbits within each are necessarily in different invariant subsets of \(\bar{E}\). White stripes indicate that the corresponding orbit is disconnected from a large part of \(\bar{E}\). Disconnected regions directly impact the mixing dynamics by containing the mixing to a region smaller than the entire exceptional set.

The larger the value of \(\rho_{ij}\) [darker in Fig.~\ref{fig:45_return_and_detail}(a)], the more often an orbit returns to that bin. If an orbit spends many iterations returning to the same cutting lines, it is not spending as many iterations exploring the domain. As such, a high value of \(\rho_{ij}\) indicates the associated orbit is small or short. Similarly, large orbits that spread throughout the domain return to the cutting lines with a much lower frequency (lighter color). We refer to filled regions in the same horizontal (row) or vertical (column) position together as \textit{blocks} due to their rectangular shape. Each block has orbits that return close to everywhere else in the block, indicating potential mixing. Since mixing can only occur between connected regions, disconnected blocks represent disparate mixing regions.

Figure~\ref{fig:45_return_and_detail}(a) shows that some regions (invariant sets) are completely disconnected from the rest of \(\bar{E}\), as indicated by the white stripes in the return plot. For example, in the lower left (corresponding to \(\mathcal{D}_2\), \(\theta \in [0,2\pi)\)) and upper right (corresponding to \(\mathcal{D}_1\), \(\theta \in [2\pi,4\pi)\)) the dark square blocks with empty vertical and horizontal stripes around them correspond directly with the blue and red arrowheads in Fig.~\ref{fig:45combined}(a). These are regions that only interact with one of the two cutting lines and are isolated from the rest of \(\bar{E}\).

The return plot in Fig.~\ref{fig:45_return_and_detail}(a) can be used to split the exceptional set into its constituent invariant sets. To explain, consider a detail of the lower left corner of Fig.~\ref{fig:45_return_and_detail}(a), which is shown in Fig.~\ref{fig:45_return_and_detail}(b) for \(\theta \in [0,\pi/4]\).
A close examination of Fig.~\ref{fig:45_return_and_detail}(b) demonstrates how the repeating section shown in the figure approximately decomposes the exceptional set for the \((45^\circ,45^\circ)\) protocol into different invariant sets. The thin curve at the center of the zigzag in the exceptional set, shown in Fig.~\ref{fig:45_return_and_detail}(e), intersects \((\theta_i,\theta_j) = (0,0)\) in the domain shown in Fig.~\ref{fig:45_return_and_detail}(b).
Immediately surrounding this zigzag curve is a caterpillar-like shell in Fig.~\ref{fig:45_return_and_detail}(h) and a region outside the arrowheads and the zigzag in Fig.~\ref{fig:45_return_and_detail}(c).
Some invariant sets only intersect the blue cutting line \(\mathcal{D}_2\). A small region inside the caterpillar-like zigzag is distinct from the central zigzag curve, as shown in Fig.~\ref{fig:45_return_and_detail}(d). The boundary of the blue arrowheads in Fig.~\ref{fig:45_return_and_detail}(g) is actually a line of small dark blocks, each representing a thin invariant set that outlines an arrowhead, as shown in Fig.~\ref{fig:45_return_and_detail}(g) surrounding the blue arrowhead core in Fig.~\ref{fig:45_return_and_detail}(f). Finally, there is a separate set intersecting \((\theta_i, \theta_j) = (\pi/4, \pi/4)\) (not shown) which is a point of tangency with a circular cell and produces a series of very small circles throughout the exceptional set. The associated short orbits that return frequently to the cutting line create the extreme skew in the return plot, which in turn requires the use of the \(\log\) scale in the figures.

None of the gray blocks associated with the \((45^\circ, 45^\circ)\) protocol in Fig.~\ref{fig:45_return_and_detail}(b) have a uniform \(\rho_{ij}\) throughout the block. The existence of nonuniform patterns within blocks indicates that these blocks are not the smallest invariant sets, although different sets may be so intertwined with one another that separating them in this fashion is impossible. This is obvious in the upper right of Fig.~\ref{fig:45_return_and_detail}(b), where the cutting line density has an almost checkerboard pattern suggesting that different parts of the arrowhead in Fig.~\ref{fig:45_return_and_detail}(f) correspond to different intensities in the block. On the other hand, the combination of two blocks is evident in the color variation within Fig.~\ref{fig:45_return_and_detail}(c), demonstrating how blocks may correspond to intertwined invariant sets that are close but not overlapping.

In order for a mixing block to exist, an orbit in the block must intersect a non-zero length of the cutting line. Combined with an occupation of a non-zero width away from the cutting line in an orbit, this means that any mixing block must have orbits that individually have non-zero area on the hemispherical shell (finite length and width). It appears that a necessary condition for mixing is that individual orbits occupy a finite area, as it allows orbits to interact with one another while thin (zero-area) curves that do not overlap anywhere do not allow mixing between them. Although a large exceptional set can be constructed using thin curves (see Appendix~\ref{sec:app_C_example_density} in the Supplementary Material), there will be no mixing if orbits do not create the area for mixing to occur.

With this decomposition of the \(\theta_i,\theta_j \in [0,\pi/4]\) region in mind, the entirety of the structure in Fig.~\ref{fig:45_return_and_detail}(a) can be decomposed into similar invariant sets. Invariant sets intersecting both \(\mathcal{D}_1\) and \(\mathcal{D}_2\) from \([0,\pi/4] \times [0,\pi/4]\) are repeated throughout the entire \([0,4\pi] \times [0,4\pi]\) space. Invariant sets intersecting just the blue cutting line, \(\mathcal{D}_2\), are repeated only in the lower left, \([0,2\pi)\times [0,2\pi)\). The invariant sets intersecting just the red cutting line, \(\mathcal{D}_1\), are repeated only in the upper right, \([2\pi,4\pi)\times [2\pi,4\pi)\), but are identical in their return structure to those that intersect just \(\mathcal{D}_2\).

\begin{figure*}  
  \includegraphics[width=0.8\textwidth]{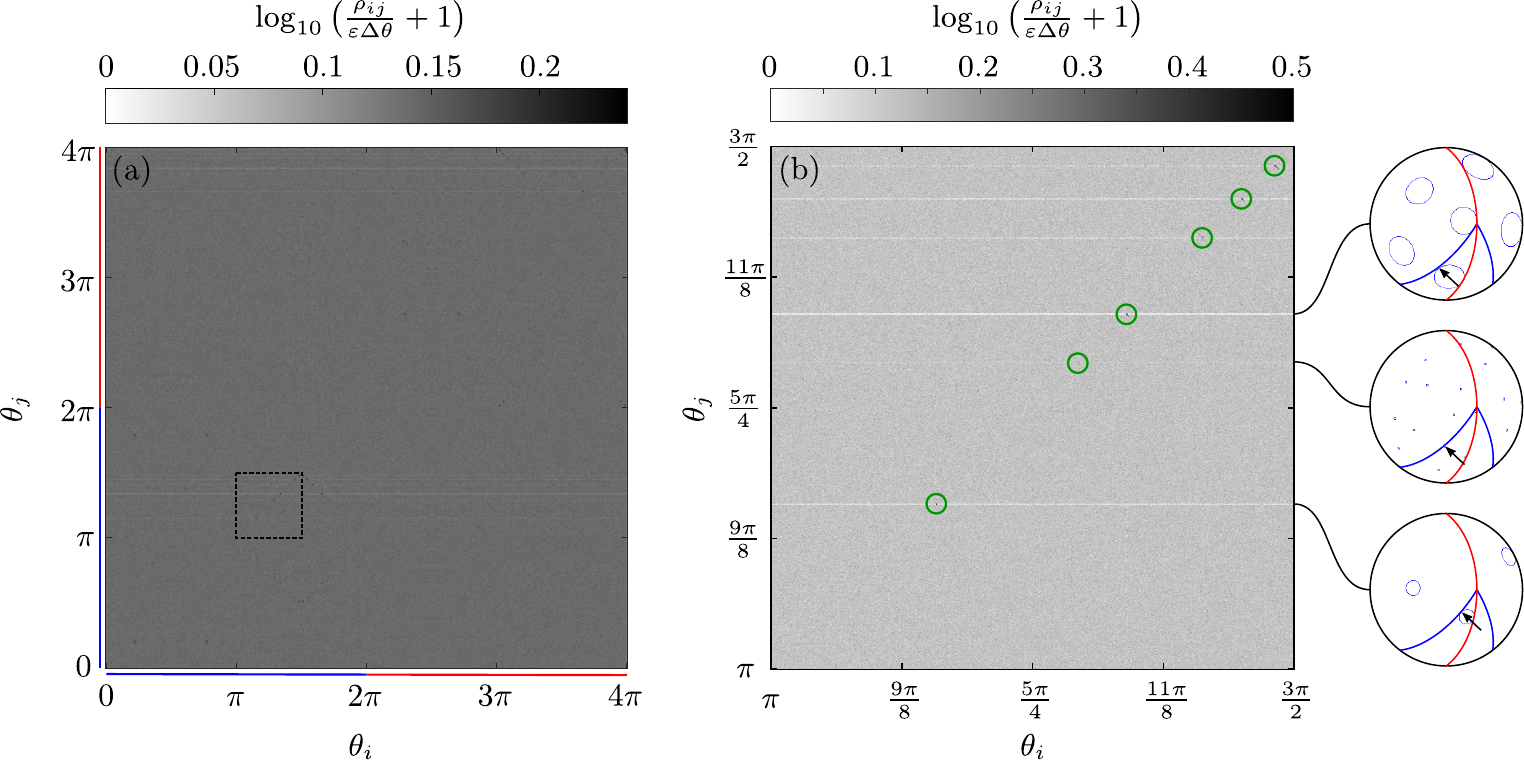}
  \caption{
  (a) Approximate return plot for the \((57^\circ, 32.75^\circ)\) protocol as a 2D histogram of \(\rho_{ij}\), return to bin \(i\) from bin \(j\); \(\Delta\theta = \pi\times 10^{-3}\). (b) Detail of (a) from the dashed black box in the lower left; \(\Delta\theta = \pi/8\times 10^{-3}\). Some of the dark spots along the diagonal are circled in green with the corresponding cells responsible for the horizontal white lines and dark spots shown on the right. The arrows indicate the point of adjacency with the circular cells responsible for the depletion of area. (a,b) \(2 \times 10^6\) points are initialized \(\delta = 10^{-9}\) from the cutting line, and iterated for \(N = 2\times 10^6\) iterations using a cutting line width \(\varepsilon = 10^{-5}\). The \(\theta \in [0,2\pi)\) range represents \(\mathcal{D}_2\) (blue/black cutting line) and the \(\theta \in [2\pi,4\pi)\) range represents \(\mathcal{D}_1\) (red cutting line).  
  }
	\label{fig:57_return_and_detail}
\end{figure*}

The return plot for \((45^\circ, 45^\circ)\) corresponds with the structure evident when invariant subsets are color-coded as in Fig.~\ref{fig:45combined}. In contrast, the return plot for \((57^\circ, 32.75^\circ)\), shown in Fig.~\ref{fig:57_return_and_detail}(a), has no distinct blocks that would indicate distinct invariant subsets. This is consistent with the previous observations that there are no distinct invariant subsets within \(\bar{E}\) in Fig.~\ref{fig:45combined}(c). In fact, if \(\bar{E}\) for \((57^\circ,32.75^\circ)\) is truly ergodic, then Fig.~\ref{fig:45combined}(c) would be uniform except for the non-mixing cells. However, upon close inspection of a small section of the return plot for the \((57^\circ, 32.75^\circ)\) protocol in Fig.~\ref{fig:57_return_and_detail}(b), faint white horizontal lines are visible which result from the proximity of orbits to a large circular periodic cell. In addition, faint dark dots, inside the superimposed circles, are also visible along each horizontal white line at the intersection with the \(\theta_i = \theta_j\) diagonal. These dark spots are from points placed \(\delta\) away from \(\mathcal{D}\) that fall inside a periodic cell. The horizontal white lines result from the decreased width of \(\bar{E}\) immediately next to a cell. Essentially, the cell occupies some of the \(\mathcal{D}_\varepsilon\) sample region which creates a light region since points are less likely to be in this small region compared with the full or nearly full width regions elsewhere along the cutting line. White vertical lines from points placed in a cell are noticeably absent since there are several hundred points seeded in each horizontal bin. Since \(\delta \ll \varepsilon\), the vertical stripes produced by placement inside a cell (due to \(\delta\)) are smaller (and not visible) than those produced by the intersection of a cell with \(\mathcal{D}_\varepsilon\) (due to \(\varepsilon\)). Thus, the trapping of points within cells is not the dominant effect producing these aberrations, but rather it is the depletion of area of \(\mathcal{D}_\varepsilon \cap \bar{E}\) due to adjacent circular cells which makes returning to cell-adjacent regions less likely.

\begin{figure}
    \includegraphics[width=0.36\textwidth]{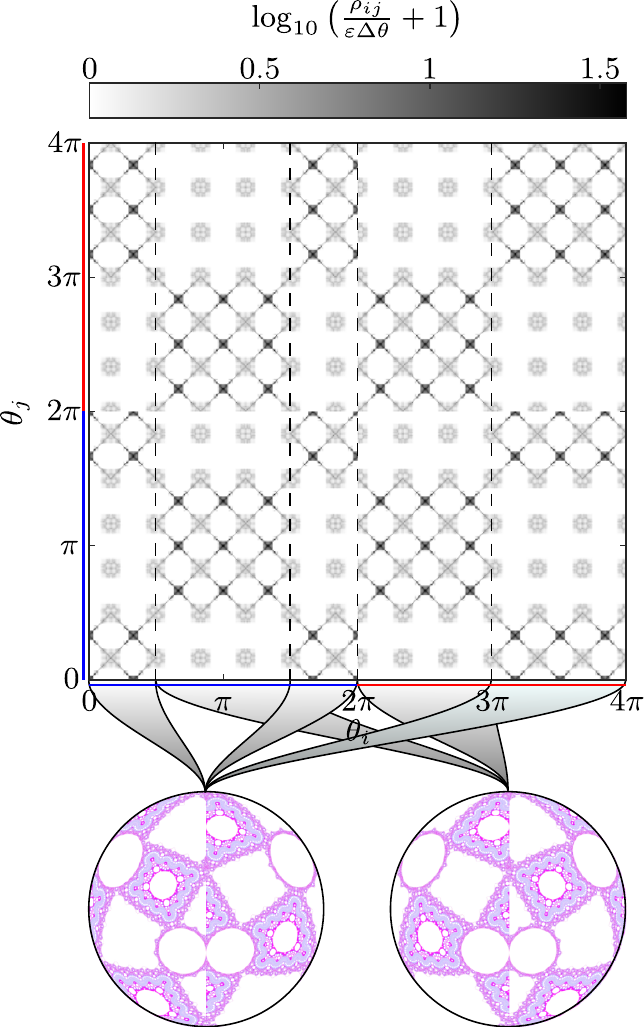}
	\caption{Return plot (\(\delta = 10^{-9}\)) for \((90^\circ, 60^\circ)\) protocol showing \(\rho_{ij}\), using \(2\times 10^6\) points, \(N = 2\times 10^6\) iterations, and \(\varepsilon = 10^{-5}\).}
	\label{fig:90-60-90_return}
\end{figure}

Recall that this `sided' parameterization of the cutting lines was chosen to specifically deal with cases where cutting lines are coincident, such as the \((90^\circ, 60^\circ)\) protocol shown in Fig.~\ref{fig:906090_PWI}. The corresponding return plot is shown in Fig.~\ref{fig:90-60-90_return}. In this case, \(\theta \in [2\pi,3\pi)\) and \(\theta \in [3\pi,4\pi)\) are two `sides' of the same cutting line, each corresponding to different atoms for the PWI dynamics. These two sides act in a similar fashion to the two separate cutting lines in the \((45^\circ, 45^\circ)\) protocol. For example, the symmetry across \(\theta = 3\pi\) for \(\theta \in [2\pi,4\pi)\) evident for the \((45^\circ, 45^\circ)\) protocol is absent in Fig.~\ref{fig:90-60-90_return}. The two distinct regions (each with the same vertical/horizontal pattern) can also be used to separate the exceptional set into two portions, as shown below the return plot. These two portions are identical when reflected about the vertical diameter, consistent with the left-right barrier to mixing shown in Fig.~\ref{fig:906090_PWI}(b). Not only is there a barrier to mixing across the vertical midline, but there are also completely separate regions contributing to each half of \(\bar{E}\) that are not captured by the red and blue coloring for the cutting lines used in Fig.~\ref{fig:906090_PWI}(c). Based on the return plot, it becomes clear that the red and blue cutting line assignment for color is naive in this case and does not capture the separable mixing sets. Separable mixing sets in \(\bar{E}\) [like those for the \((45^\circ, 45^\circ)\) protocol in Fig.~\ref{fig:45_return_and_detail}(b)] are not a generic property. Instead, they occur for some protocols but not for others. However, separable mixing sets are a generic property of protocols with exclusively coincident cutting lines (polygonal tilings) as shown in Appendix~\ref{sec:app_D_iet_protocols} in the Supplementary Material.

\begin{figure*}
  \includegraphics[width=0.95\textwidth]{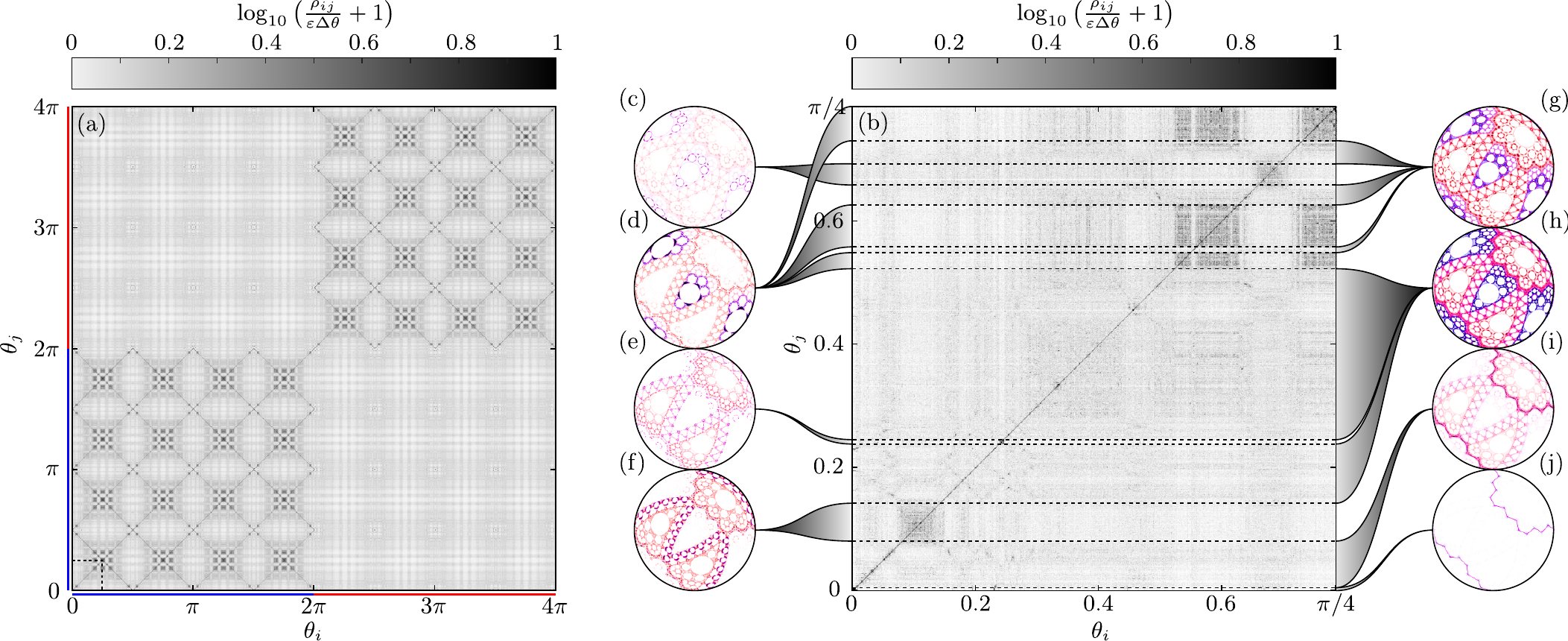}
  \caption{
  (a) Approximate return plot for the \((57^\circ, 57^\circ)\) protocol as a 2D histogram of \(\rho_{ij}\), return to bin \(i\) from bin \(j\); \(\Delta\theta = \pi\times 10^{-3}\). (b) Detail of (a) from the dashed black box in the lower left; \(\Delta\theta = \pi/16\times 10^{-3}\). (a,b) \(2 \times 10^6\) points are initialized \(\delta = 10^{-9}\) from the cutting line, iterated \(N = 2\times 10^6\) times using a cutting line width \(\varepsilon = 10^{-5}\). (c-j) Various regions of the exceptional set corresponding to portions of the return plot.
  }
  \label{fig:return_57-57}
\end{figure*}

We have described an easily separable protocol that has large regions with no returns at all (white) in the return plot indicating barriers to mixing, \((45^\circ, 45^\circ)\), and a seemingly inseparable protocol that has no large regions without returns indicating few barriers to mixing, \((57^\circ, 32.75^\circ)\). An intermediate case is the \((57^\circ,57^\circ)\) protocol, first shown in Fig.~\ref{fig:pwi_demo}(g), for which only 0.4\% of the return plot in Fig.~\ref{fig:return_57-57}(a) is empty (white), but which still has many barriers to mixing due to separate invariant subsets within the exceptional set. This is made apparent by the obvious pattern in this protocol's return plot in Fig.~\ref{fig:return_57-57}(a). The comparison between the \((57^\circ, 57^\circ)\) protocol and the other two we have examined is imperfect because the fractional coverage for the \((57^\circ, 57^\circ)\) protocol is \(\Phi_\infty \approx 0.33\), smaller than the coverage of the other two. Nevertheless, the \((57^\circ, 57^\circ)\) protocol has orbits that are not equally dense in their returns to the cutting lines, indicating that orbits are simultaneously very close to one another and still separate.

To illustrate, the \([0,\pi/4]\times[0,\pi/4]\) region of Fig.~\ref{fig:return_57-57}(a) is shown in detail and approximately decomposed into different invariant subsets in Fig.~\ref{fig:return_57-57}(b). Although there is a pattern within the return plot, it is essentially blurred or smudged such that the more precise separation done in Fig.~\ref{fig:45_return_and_detail}(b) for the \((45^\circ, 45^\circ)\) protocol is not possible here. Instead, some of the more `fuzzy' regions can be approximately separated. Figures~\ref{fig:return_57-57}(c,d) approximate the invariant sets that compose the center of the arrowhead structures visible in Fig.~\ref{fig:pwi_demo}(g). That the colors are very faint is a consequence of almost all orbits reaching almost all of the exceptional set with a low density even when each orbit has a region of high density. In order to better show these more disperse, faint orbits, the colors used for the orbits in Fig.~\ref{fig:return_57-57}(c-j) are twice as dark as those used in Fig.~\ref{fig:45_return_and_detail}(b). Figures~\ref{fig:return_57-57}(e) and (f) show regions of slightly concentrated blue and wide faint regions of red, exploring the boundary between the major blue and red arrowheads. The regions in (g) and (h) explore most of the exceptional set uniformly other than a difference in blue and red density. Figure~\ref{fig:return_57-57}(i) shows a region surrounding the zigzag curve shown in (j), similar to the zigzag in Fig.~\ref{fig:45_return_and_detail}(e) for the \((45^\circ, 45^\circ)\) protocol. Thus, the \((57^\circ, 57^\circ)\) protocol return plot represents an intermediate between the easily separable \((45^\circ, 45^\circ)\) protocol and the seemingly impossible-to-separate \((57^\circ, 32.75^\circ)\) protocol in which almost all orbits visit most of the exceptional set with uneven density. 

The return plot contains both information about the area and connectivity of the potentially ergodic region of the PWI, \(\bar{E}\). In an abstract sense, the return plot defines the connectivity of a complex network between different segments of the cutting lines. This return plot analysis is naturally extendable to any area preserving PWI with finite perimeter atoms (to allow parameterization) and could have applications in the analysis of other area preserving dynamical systems. We conjecture that a metric on the shapes within the return plot in conjunction with the fractional coverage \(\Phi_\infty\) would create a complete measure of mixing for any given PWI of this type. In the following section, we show that \(\Phi_\infty\) can be extracted from density measurements such as the return plot. Variation within the return plot seems to be an indicator of distinct mixing regions that form barriers to overall mixing. The amount of white space in the return plot also is an indicator of the lack of connectivity of the exceptional set that indicates the barriers to mixing within \(\bar{E}\). A measurement of the variation in density or the amount of white space (disconnectedness) could, in conjunction with a measurement of the area of \(\bar{E}\), be used to quantify the mixing under any given protocol.

\section{Estimating mixing metrics from cutting line density values} \label{sec:density_integrals}

We now return to the fractional coverage of the fat fractal exceptional set, \(\Phi_\infty\) \cite{Park2017,Smith2017}, which is the normalized area of the fat fractal, to show that the density measure described in Section~\ref{sec:natural_invariant} can be used to estimate the area of the fat fractal, and hence the value of \(\Phi_\infty\), from measurements of \(\rho\) along the cutting lines. One of the more interesting features of the exceptional set is that a two-dimensional fat fractal is generated from one-dimensional  cutting lines. Since the PWI is composed of isometries, this means that all information (without distortion) about the exceptional set, specifically the cutting line density everywhere, must be contained within the original one-dimensional, \(\varepsilon\)-fattened cutting lines themselves, and applying the PWI over many iterations simply spreads this information to the rest of the exceptional set. With this in mind, the cutting line density \(\rho(\bm{x},\mathcal{D}_\varepsilon)\) at any point along \(\mathcal{D}_1\) or \(\mathcal{D}_2\) is a measure of how that particular point on the cutting line spreads through the domain under the PWI.

Since the PWI map is measure preserving, \(\rho(\bm{x},A)\) can be extended from a simple counting measure fraction to a fraction of how much of the orbit \(X = \bigcup_{i = -\infty}^\infty M^i (\bm{x})\) is in an area \(A\), i.e.\ a fraction of measures \(\rho(\bm{x},A) = \mu(A\cap X)/\mu(X)\). The Lebesgue measure (area) of an orbit \(X\) from point \(\bm{x} \in \mathcal{D}_\varepsilon\), where \(\mathcal{D}\) is any one-dimensional curve such as the cutting lines (or any subset of the cutting lines), is
\begin{equation}
\mathcal{L}^2(X) = 
	\int_{\mathcal{D}_\varepsilon \cap X}  \frac{ 1}{\rho(\bm{x},\mathcal{D}_\varepsilon)} d\bm{x}.
	\label{eq:single_orbit_size_body}
\end{equation}
This integral is specifically over a two dimensional region such that if the dimension of the orbit is smaller than two, this integral evaluates to zero. Since \(\rho(\bm{x},\mathcal{D}_\varepsilon)\) is invariant under the map, \(1/\rho(\bm{x},\mathcal{D}_\varepsilon)\) is also invariant. Intuitively, the integral bounds measure the size of the intersection of \(X\) and \(\mathcal{D}_\varepsilon\) and the value of \(1/\rho(\bm{x},\mathcal{D}_\varepsilon)\) measures how many copies of this intersection are used to create the orbit \(X\).

By taking the union of all orbits from \(\mathcal{D}_\varepsilon\), which is precisely the \(\varepsilon\)-fattened exceptional set \(E_{\varepsilon}\), the area of the \(\varepsilon\)-fattened exceptional set is
\begin{equation}
  \mathcal{L}^2(E_\varepsilon) =  \int_{\mathcal{D}_\varepsilon}  \frac{1}{\rho(\bm{x},\mathcal{D}_\varepsilon)} d\bm{x}. \label{eq:exceptional_set_size_body_exact}
\end{equation}
The intersection of \(\mathcal{D}\) and any individual orbit \(X\) is no longer needed as the intersection of all orbits from \(\mathcal{D}_\varepsilon\) with itself is simply \(\mathcal{D}_\varepsilon\). In this way, the integral evaluates the expected value of the number of copies of \(\mathcal{D}_\varepsilon\) needed to create \(E_\varepsilon\) and multiplies this by the area of \(\mathcal{D}_\varepsilon\) to get the total area.

In the limit of small \(\varepsilon\), the region of integration approaches \(\mathcal{D}\) which is one dimensional. With this in mind, the integral in Eq.~\ref{eq:exceptional_set_size_body_exact} can be approximated by a one dimensional integral over \(\mathcal{D}\) multiplied by width \(\varepsilon\),
\begin{equation}
  \mathcal{L}^2(\bar{E}) \approx \lim_{\varepsilon \to 0} \int_{\mathcal{D}_a}  \frac{ \varepsilon}{\rho(\bm{x},\mathcal{D}_\varepsilon)} d\bm{x} + \lim_{\varepsilon \to 0} \int_{\mathcal{D}_b}  \frac{ \varepsilon}{\rho(\bm{x},\mathcal{D}_\varepsilon)} d\bm{x}. \label{eq:exceptional_set_size_body}
\end{equation}
where \(\mathcal{D}_a\) and \(\mathcal{D}_b\) are the two \textit{sides} of the cutting line as shown in Fig.~\ref{fig:cutting_line_assembly}, each side is assigned the dynamics of one of the two overlapping atoms, \(P_a\) or \(P_b\). The value of  \(\lim_{\varepsilon \to 0} \varepsilon/\rho(\bm{x},\mathcal{D}_\varepsilon)\) is an abstraction of how much of the orbit is spread away from the measurement set. Note that for periodic orbits, \(\lim_{\varepsilon \to 0} \rho(\bm{x},\mathcal{D}_\varepsilon) \neq 0\) such that \(\lim_{\varepsilon \to 0} \varepsilon/\rho(\bm{x},\mathcal{D}_\varepsilon)\) is always 0. To arrive at \(\Phi_\infty\), this measure simply needs to be normalized by the area of domain. Since the \(\varepsilon \to 0\) limit is not directly computable, an appropriately small \(\varepsilon\) can be used to approximate \(\mathcal{L}^2(\bar{E})\). More thorough derivations of Eqs.~\ref{eq:single_orbit_size_body}-\ref{eq:exceptional_set_size_body} as well as additional intuition regarding density and area are given in Appendix~\ref{sec:app_B_integral_proof} in the Supplementary Material.

The parameterization in Fig.~\ref{fig:parameterization} provides a framework for applying the integral expression for the area of the exceptional set given in Eq.~\ref{eq:exceptional_set_size_body}, such that this integral becomes
\begin{equation}
	\mathcal{L}^2(\bar{E}) \approx \lim_{\varepsilon \to 0} \int_{0}^{4\pi} \frac{\varepsilon}{\rho(\bm{s}(\theta),\mathcal{D}_\varepsilon)}\ d\theta. \label{eq:final_integral_paramet}
\end{equation}
The previously used metric for the fractional coverage of the exceptional set was \(\Phi_\infty = \mathcal{L}^2(\bar{E}) / 2\pi\), where \(2\pi\) is the area of the unit hemispherical shell \cite{Park2017, Lynn2019, Smith2017}. The density along the cutting lines \(\rho(\bm{s}(\theta),\mathcal{D}_\varepsilon)\) is, in general, not continuous, so it is convenient to recast the integral as an average (denoted with angled brackets) of the integrand over \(\theta\) to evaluate \(\Phi_\infty\),
\begin{eqnarray}
	\Phi_\infty &\approx& \frac{1}{2\pi} \lim_{\varepsilon \to 0} \int_{0}^{4\pi} \frac{\varepsilon}{\rho(\bm{s}(\theta),\mathcal{D}_\varepsilon)}\ d\theta \\
	&\approx& \lim_{\varepsilon \to 0} \left\langle \frac{2 \varepsilon}{\rho(\bm{s}(\theta),\mathcal{D}_\varepsilon)} \right\rangle. \label{eq:average}
\end{eqnarray}

From the return plots in Sec.~\ref{sec:return_maps}, \(\rho(\bm{s}(\theta),\mathcal{D}_\varepsilon)\) can be approximated by integrating (summing) any corresponding column (or row) in the return plot, which reflects the total return to (or from) any part of the cutting line, i.e.\ summing across row \(j\) approximates \(\rho(s(\theta \in [\theta_i,\theta_i+\Delta \theta)),\mathcal{D_\varepsilon}) \approx \sum_{j} \rho_{ij} \) and the fractional coverage would be \(\Phi_\infty \approx \lim_{\varepsilon\to 0}\langle 2\varepsilon (\sum_i \rho_{ij})^{-1}\rangle_j\).
For instance, the \((57^\circ,32.75^\circ)\) protocol has a value of \(\Phi_\infty \approx 0.41\) using direct measurement \cite{Lynn2019} and \(\Phi_\infty \approx 0.41\) using Eq.~\ref{eq:average} with \(\varepsilon = 10^{-5}\), \(\delta = 10^{-9}\), \(N = 2\times 10^6\) iterations, averaged across \(2\times 10^6\) different positions along the cutting lines.

\section{Conclusions}
We have described two different representations of the internal dynamics in the exceptional set of the hemispherical PWI by creating an invariant measure, \(\rho\), which we refer to as \textit{cutting line density}. The first representation is easily shown everywhere throughout the domain through a color-coding that is generated by interactions with the two different cutting line discontinuities in the PWI. The second is the natural extension of this color-coding, treating the cutting lines as a continuum, and measuring the return dynamics at each point along them to create a \textit{return plot}. Although the color-coding can show a snapshot of the dynamics everywhere in the exceptional set, the return plot only shows dynamics near the original cutting lines. Due to the nature of the construction of the PWI exceptional set, the dynamics at the cutting line are representative of the dynamics everywhere in the exceptional set.

These two representations of the cutting line density, color-coding and return plots, indicate that more cuts do not indicate more mixing, rather an even distribution of cuts is required for mixing and variation in cutting line densities actually indicate barriers to mixing. Additionally, sufficient variation in the return dynamics of the cutting line are required to generate mixing. Without this variation, no amount of cuts can produce mixing, even if the domain is cut with an equal density of cuts. Ultimately, the existence and detection of confined mixing regions in the spherical tumbler mixer is aided by understanding the role of cutting-and-shuffling in creating trapped orbits. Large, invariant regions in the PWI, such as the arrowheads in the \((45^\circ, 45^\circ)\) protocol, as shown in Fig.~\ref{fig:coverage_compare}(a), are essentially extensions of the non-mixing cells contained in them.

The density of cutting lines in the hemispherical shell PWI can be used as an indicator of non-ergodic behavior within the exceptional set. Although the fractional coverage of the exceptional set, \(\Phi_\infty\), is correlated with the degree of mixing \cite{Park2017,Smith2017}, the existence of mixing within the exceptional set is entirely dependent on the dynamics within it. The cutting line density based coloring techniques described in Section \ref{sec:coloring} can be used in most situations to rule out (but not to definitively show) ergodicity within an exceptional set with minimal effort.

Return plots described in Section \ref{sec:return_maps} essentially project 2D orbits into a 1D space (the parameterized cutting lines), which necessarily exists for this PWI since the 1D set of cutting lines \(\mathcal{D}\) sweeps out and generates the exceptional set. Separating out orbits using the return plot can reveal non-mixing regions within the exceptional set as well as rule out ergodicity. A measurement of the white space (disconnectedness) of the return plot or the variation of the local densities within the plot, combined with the fractional coverage of the exceptional set, indicates not only how much of the domain can mix, but also accounts for how mixing within this space is occurring. Return plots can also be used to determine cutting line densities from which \(\Phi_\infty\) can be approximated. Consequently, return plots enable relevant mixing information to be contained in a single representation that expresses both the size and connectivity of the exceptional set.

Just as cells only interchange with other cells within their periodic orbit \cite{Goetz2003, Smith2017, Park2016a, Scott2003, Bruin2003}, inside the exceptional set points only mix with one another if they are connected by at least one orbit. The colored exceptional set and the return plot are both visual representations of the connectedness of all the orbits within the exceptional set. The color-based classification of the exceptional set allows easy identification of invariant subsets, while the return plot provides a detailed deconstruction of where these subsets appear in the exceptional set. The return plot image is similar to that of an adjacency matrix in a weighted graph system if such a graph had infinitely many nodes \cite{Ray2013}.

The methods presented in this paper to examine mixing, connectedness, and barriers to both mixing and connectedness can be immediately applied to any 2D, finite domain, almost everywhere invertible PWI (especially different geometries \cite{Smith2019}) and even PWI systems of other dimensions (although a 2D return plot may not always be possible). Notably, these methods could be easily applied to simpler PWIs such as a digital filter (a square geometry) \cite{Ashwin1997a, Ashwin1997}. Apart from the change in geometry, little changes in the analysis. The ability to construct and measure the connectivity of the exceptional set for any such system provides an effective means of evaluating the mixing, and the coloring methods can be applied to immediately identify the existence of invariant subsets by inspection.

\section{Acknowledgments}
We gratefully acknowledge Lachlan D.~Smith and Christian Gorski for insightful discussions. This research was based upon work supported by the National Science Foundation under grant No.~CMMI-1435065.

\appendix


%


\clearpage

\section{Discussion of error}\label{sec:app_A_errors}

There are four main sources of error in the methods presented here: error due to numerical precision, error due to \(\delta\)-placement away from cutting lines (for the return plot and measurements of \(\rho\)), error due to \(\varepsilon\)-fattened cutting lines, and error due to the number of iterations, \(N\), being finite.

Numerical precision affects all numerically advanced map such as the PWI studied here. Error from finite precision calculations manifest as accumulated errors in position, or drift, that scales like \(N\) times the machine epsilon (in this case, roughly \(2\times 10^{-16}\)). Better precision in position is not possible using double-precision data-types since this error is scaled by the position of the tracked point, which is often \(\mathcal{O}(1)\). These types of numerical errors could influence the underlying dynamics of the map by creating spurious new periodic points in the numerical system or destroying actual periodic orbits through drift. The scale of these types of errors is mostly irrelevant to the study here (other error effects are much larger in comparison) because the accumulated drift after \(2\times 10^6\) iterations is only about \(4\times 10^{-10}\).

\(\delta\)-placement away from the cutting lines is introduced to ensure that the aforementioned numerical drift does not cause a calculated point to incorrectly drift across the cutting line it is placed near. The value in this paper, \(\delta = 10^{-9}\), is chosen to account for an \(N\times (\text{machine epsilon}) = 4\times 10^{-10}\) amount of drift. The numerical drift can still cause errant cutting line crossings at later positions, but a positive \(\delta\) guarantees that points do not immediately collapse onto one another (removing the sided nature of the cutting lines). If a point is advected under the map and then advected under the reverse map, this \(\delta\)-placement should keep the original point and the iterated point in the same starting atom. The biggest consequence of \(\delta\)-placement is the potential placement of points inside of a cell instead of inside the exceptional set. The fraction of points errantly placed inside of a cell is incalculable (it depends on the details of \(\bar{E}\)), but points placed within cells can only occupy thin circular orbits. These thin circular orbits contribute a \(\delta\)-fattening to the exceptional set when taken together. In addition to this fattening, they can also exclude parts of the exceptional set; the exceptional set is only guaranteed to be arbitrarily close to the cutting lines, not necessarily arbitrarily close to a line \(\delta\) from the cutting lines. Any invariant sets within the exceptional set not intersecting this line are not measured. Hence, the error from \(\delta\) changes values of \(\rho(\bm{x},A)\) by shifting \(\bm{x}\) by \(\delta\). The balance of unmeasured invariant sets and a measured fattening is difficult to determine. However, the error due to \(\delta\)-placement is potentially orders of magnitude smaller than the error incurred from \(\varepsilon\)-fattening of the cutting line.

The error due to \(\varepsilon\)-fattening introduces an additional, adjacent region, where points can `return' to the cutting line, and creates a border within the collection of periodic islands, \(O\), where points in open cells appear to belong to the exceptional set. In the colored exceptional set [e.g., Fig.~\ref{fig:pwi_demo}(g)], these regions appear as distinct invariant sets, which confuses any measurement of ergodicity (which would only have a single invariant set within \(\bar{E}\)). In the return plot, this \(\varepsilon\)-fattening causes a vertical smearing effect since orbits can return to other parts of the cutting lines without actually reaching them. This smearing is not uniform, as evidenced by the horizontal white bands, discussed with respect to Fig.~\ref{fig:57_return_and_detail}, which result from the overlap with the circular cells in \(O\). For direct measurements of \(\Phi\), the fattening of the cutting lines includes the entirety of any cells of radius \(\le \varepsilon\) and the boundary of all cells with larger radii. For the density based calculation of \(\Phi\), error is not primarily dependent on \(\varepsilon\) but rather on the placement of points \(\delta\) from the cutting line and the spacing between points (an accurate measurement would require values of \(\rho\) continuously along the cutting line). As outlined in Appendix~\ref{sec:app_B_integral_proof}, the size and shape of the region used for measurement of \(\rho\) is not as important as the intersection with invariant sets, and the \(\delta\)-placed line is solely responsible for this. However, there is error in measuring the size of the fattened cutting line (since the line lies on a sphere, the cutting line is a great circle, but the \(\varepsilon\)-fattened border is not) and overlap between the two cutting lines. This error can be quantified exactly, although the exact overlap is not easily formulated \footnote{The exact value of this overlap is \(4\int_0^{\frac{\pi}{2}}4-4\sqrt{1-\varepsilon^2 \csc(\beta + \phi)^2} -4\sqrt{1-\varepsilon^2 \csc(\alpha + \phi)^2}d\phi\) when \(\varepsilon \le \min[\sin(\alpha),\sin(\beta)]\).}. The approximate area of the fattened cutting line is \(4\pi\varepsilon\) while the true area is \(4\pi\sin(\varepsilon)\) minus the overlap of the two fattened cutting lines which is asymptotically \(2\varepsilon^2[\csc(\alpha)\sec(\alpha) + \csc(\beta)\sec(\beta)]\) for small \(\varepsilon\). The overlap is the dominant error source since \(\varepsilon - \sin(\varepsilon) = \mathcal{O}(\varepsilon^3)\) and the overlap is \(\mathcal{O}[\varepsilon^2(\frac{1}{\alpha} + \frac{1}{\beta})]\) accounting for the \(\csc\) terms which are larger when \(\alpha\) or \(\beta\) are small. The error due to \(\varepsilon\) changes the value of \(\rho(\bm{x},A)\) by introducing overlap at the cutting line intersection (double-counting) but is still small, \(\mathcal{O}(\varepsilon^2)\), for most of the \((\ab)\) protocol space.

Error due to finite iterations is subtle since it depends sensitively on the irrationality of measured quantities. Irrationality here concerns rotations and measurements that are not `close' to rational numbers that are `early' in Cantor's ordering of the rationals. For example, a cell that rotates internally by an amount that is nearly \(1/3\) will appear as a triangle for much longer than a cell that rotates internally by an amount that is nearly \(1/300\), hence the \(1/300\) rotation may appear more irrational since more of the boundary of the cell is initially explored. Quantifying this irrationality is subjective and dependent on the values of \(\varepsilon\) and \(N\). If, in the previous example, the rotation is orders of magnitude closer to \(1/300\) for the second cell than \(1/3\) of the first, which one appears more irrational may switch after some \(N\). These are errors in the construction of \(\tilde{E}_{\varepsilon,N}\), but there are also errors in measuring \(\rho\). If \(\rho\) at some point is truly irrational, then there is no sufficient \(N\) at which \(\rho\) can be measured exactly. In this case, the true value of \(\rho\) can be considered the average of a random variable, i.e. \(\rho(\bm{x},A)\) is the probability of the orbit of \(\bm{x}\) being in \(A\), for which one can consider the standard error of the mean, which is \(\sqrt{\rho(1-\rho)/N}\). The \(\sqrt{N}\) dependence would be a dominant effect if not for the influence of the actual value of \(\rho\), which is small for non-periodic orbits and small \(\varepsilon\). The error due to finite \(N\) changes the value of \(\rho(\bm{x},A)\) by introducing error in the actual measurement of \(\rho\) rather than an error in orbit (\(\delta\)) or orbit intersection with the overlapped region (\(\varepsilon\)). Additionally, the error in measurements of \(\rho\) as a result of using a finite number of iterations manifests itself as a sparse return plot, and, as \(\varepsilon\) is decreased, the number of iterations required to offset this sparsity increases.

With all this said, the dominant error term is unclear since the various non-numerical errors depend on the geometry of the exceptional set which varies continuously through the \((\ab)\) protocol space. Nevertheless, our experience in considering a range of protocols with the \(N\), \(\delta\), and \(\varepsilon\) values used here, gives us confidence that the errors in the current analysis are inconsequential to the conclusions about invariant mixing regions within the exceptional set.

\section{Extracting the measure of the exceptional set from cutting line density measurements}\label{sec:app_B_integral_proof}

Previous papers \cite{Park2017, Lynn2019, Smith2017} have examined the correlation between mixing and the size of the exceptional set (the fat fractal where mixing is possible). It can be shown that the size of the exceptional set is related to the measured density of cutting lines along the initial cutting lines. Intuitively, density can only increase if cutting lines are not spreading to other parts of the hemisphere, but instead stack on themselves such that \(1/\rho\) is, abstractly, a measure of how cutting line `area' is spread through the hemisphere. Since density is easy to measure for a single point in a set, using density to approximate the measure of the exceptional set is useful.

For a measure preserving map \(M:S \to S\), define the density of the orbit from \(\bm{x}\) in set \(A\) as in Eq.~\ref{eq:density}, and define the closure of a set \(C\) as \(\bar{C} =  \lim_{r \to 0} B_r(C)\) such that the following are equivalent,
\begin{align}
\rho(\bm{x} , A) =& \lim_{r \to 0} F(\bm{x}, B_r(A)) \tag{\ref{eq:density}}\\
	=& F(\bm{x}, \bar{A}) \nonumber \\
	=& \lim_{n \to \infty} \frac{\# \{ M^i(\bm{x}) \in \bar{A} : 0 \le i \le n \}}{n} \nonumber\\
    =& \lim_{n\to \infty}\frac{\#\{M^i(\bm{x})\in \bar{A} : -n \le i < n\}}{2n},
\end{align}
where \(\#\) is the counting measure and \(\rho\) measures the fraction of iterates that an orbit spends in \(\bar{A}\). The forward iterates are sufficient for the map in this paper, but the backward iterates can be included as shown previously \cite{Lynn2019, Park2017}. To be explicit, the measurable space here is \((S,\mathcal{E})\) where \(S\) is the domain and \(\mathcal{E}\) is the invariant \(\sigma\)-algebra with respect to \(M\). Then \(\rho\) is the conditional probability of \(A\) with respect to the invariant \(\sigma\)-algebra (i.e.\ the possible combinations of minimally invariant sets under the map). That is, if \(\bm{x}\in X\) where \(X\) is a minimally invariant set, \(X \in \mathcal{E}\), such as an orbit of the map, \(X = \bigcup_{n = -\infty}^\infty M^n(\bm{x})\),
\begin{equation}
\rho(\bm{x},A) = P_{\bm{x}}(A|X) = \frac{P_{\bm{x}}(A\cap X)}{P_{\bm{x}}(X)},
\end{equation}
where \(P_{\bm{x}}(A|X)\) is the conditional probability measure such that given a random iterate of \(\bm{x}\) under the map, \(P_{\bm{x}}(A|X)\) is the probability that this iterate is in \(A\). This probability measure only makes sense if \(P_{\bm{x}}(X) \neq 0\), so we define the symbol \(P_{\bm{x}}(\cdot)\) with respect to \(\bm{x}\) such that \(P_{\bm{x}}(X)\), given \(X\) is the orbit of \(\bm{x}\), is always 1. In an inexact way, we extend this definition to measures other than \(P_{\bm{x}}(\cdot)\) and allow \(\rho(\bm{x},A) = \lim_{r \to 0} \mu(B_r(A\cap X))/\mu(B_r(X))\) where \(\mu(\cdot)\) is a higher dimensional measure matching the dimension of the ball \(B_r(\cdot)\) such that \(\mu(B_r(\cdot))\) is not zero or \(\infty\). We ultimately utilize the definition \(\rho(\bm{x},A) = \mu(A \cap X)/\mu(X)\) for intermediate steps, although this is inexact.

This density function is constant along orbits, i.e.\ it is measurable with respect to the invariant \(\sigma\)-algebra. Take \(\bm{y} \in A\) in minimally invariant set \(Y\) such that \(\rho(\bm{y},A) \neq 0,\) and therefore, since \(\rho\) is constant throughout orbits,
\begin{eqnarray}
\mu(Y) &=& \frac{\mu(A \cap Y)}{\rho(\bm{y},A)}\\
 &=& \frac{1}{\rho(\bm{y},A)} \int_{A \cap Y}d\mu \\
 &=&  \int_{A \cap Y} \frac{1}{\rho(\bm{y},A)} d\bm{y}. \label{eq:single_orbit_size_app}
\end{eqnarray}
Equation~\ref{eq:single_orbit_size_app} is the generalization of Eq.~\ref{eq:single_orbit_size_body}. Define the \textit{streak} of set \(A\), \(\mathscr{A}\), under the map as the union of all orbits to or from \(A\) and that intersect \(A\) at least once (i.e.\ the orbit of the set \(A\)), \(\mathcal{E}_i\), such that
\begin{equation}
\mathscr{A} = \bigcup_{\bm{x} \in A} \bigcup_{n = -\infty}^\infty M^n(\bm{x}) = \bigcup_{\bm{x} \in A} \mathcal{E}_i.
\end{equation}
\(\mathscr{A}\) is likewise the union of all minimally invariant sets, \(\mathcal{E}_i\), from the invariant \(\sigma\)-algebra, \(\mathcal{E}\), that intersect \(A\), i.e.\ every \(\mathcal{E}_i \cap A \neq \emptyset\) for all \(\mathcal{E}_i \in \mathcal{E}\). Likewise, \(\mathscr{A}\) is the smallest collection of minimally invariant sets \(\bigcup \mathcal{E}_i\) such that \(A \subset \bigcup \mathcal{E}_i\). To be explicit, we choose minimally invariant sets from the invariant \(\sigma\)-algebra, \(\mathcal{E}\), that intersect \(A\) and label these \(\mathcal{E}_i\). There may be uncountably many minimally invariant sets, \(\mathcal{E}_i\), that compose \(\mathscr{A}\). Then, \(\mu(\mathscr{A})\) can be written as,
\begin{equation}
\mu(\mathscr{A}) = \sum_i \mu(\mathcal{E}_i).
\end{equation}
\(\mu(\mathscr{A})\) is the `summation' of \(\mu(\mathcal{E}_i)\) for all minimally invariant sets \(\mathcal{E}_i\) that make up \(\mathscr{A}\). Since \(\bigcup_i A \cap \mathcal{E}_i = A\) and all \(\mathcal{E}_i\) are mutually disjoint, this implies
\begin{eqnarray}
\mu(\mathscr{A}) &=& \sum_i \mu(\mathcal{E}_i) \\
&=& \sum_i \int_{A \cap \mathcal{E}_i} \frac{1}{\rho(\bm{y},A)} d\bm{y} \\
&=& \int_{A} \frac{1}{\rho(\bm{y},A)} d\bm{y}. \label{eq:exceptional_set_size_app_exact}
\end{eqnarray}
Equation~\ref{eq:exceptional_set_size_app_exact} is the generalization of Eq.~\ref{eq:exceptional_set_size_body_exact}. This intuitively says that the density information within \(A\) is the reciprocal of how \(A\) is advected through the domain. If density is `high' at some \(\bm{x}\), the associated orbit is mostly contained within \(A\), while if density is `low' at some \(\bm{x}\), the associated orbit is mostly outside of \(A\), spreading out the streak of \(A\).

By referring to \(\rho\) as a density, we acknowledge that it is analogous to a physical mass density. Consider that \(A\) contributes a certain amount of `mass' to \(\mathscr{A}\) with each iteration. This contribution is constant since the map \(M\) is measure preserving (all iterates of \(A\) have the same `mass'). This means that each individual orbit from a point in \(A\) is given the same amount of mass. The density is then the distribution of this mass, and, since density is necessarily constant along orbits, the area of each orbit can be calculated, i.e.\ \(\rho = \text{mass}/\text{area}\). Note that although some orbits are given more mass with each iteration than others depending on how much the orbit overlaps the cutting line, this additional mass is taken into account by a proportional increase in \(\rho\).

This can be applied directly to fattened cutting lines \(\mathcal{D}_\varepsilon\) such that,
\begin{equation}
\rho(\bm{x},\mathcal{D}_\varepsilon) = \frac{\mu(\mathcal{D}_\varepsilon \cap X)}{\mu(X)}.
\end{equation}
Note that the streak of \(\mathcal{D}_\varepsilon\) is, for infinite iterations, the fattened exceptional set \(E_\varepsilon\), and, for finitely many iteration \(N\), this streak is the approximation to \(E_\varepsilon\), \(\tilde{E}_{\varepsilon,N}\). Note that \(\lim_{\varepsilon \to 0} (\lim_{N\to \infty} \tilde{E}_{\varepsilon, N}) = \lim_{\varepsilon \to 0} E_{\varepsilon} = \bar{E}\) \cite{Lynn2019}. Using the two dimensional Lebesgue measure instead of \(\mu(\cdot)\), the 2D measure of the approximate exceptional set is
\begin{equation}
\mathcal{L}^2(\tilde{E}_{\varepsilon,N}) \approx \mathcal{L}^2(E_\varepsilon) = \int_{\mathcal{D}_\varepsilon} \frac{1}{\rho(\bm{x},\mathcal{D}_\varepsilon)} d\bm{x}.
\end{equation}
In practice, measuring \(\rho(\bm{x},A)\) cannot be done exactly, such that this integral is closer to value of \(\tilde{E}_{\varepsilon,N}\) than \(E_\varepsilon\).

In the limit of small \(\varepsilon\), the fattened cutting line approaches a thin, one dimensional, cutting line such that \(\mathcal{L}^2(\mathcal{D}_\varepsilon) \sim 2\varepsilon \mathcal{L}^1(\mathcal{D})\) (\(\varepsilon\) contribution from either side of the cutting line). Using this, the integral is, as \(\varepsilon \to 0\), asymptotically,
\begin{eqnarray}
\lim_{\varepsilon \to 0} \mathcal{L}^2(E_{\varepsilon}) &=& \lim_{\varepsilon \to 0}  \int_{\mathcal{D}_\varepsilon} \frac{1}{\rho(\bm{x},\mathcal{D}_\varepsilon)} d\bm{x} \\
&\approx& \lim_{\varepsilon \to 0} \int_{\mathcal{D}} \frac{2\varepsilon}{\rho(\bm{x}, \mathcal{D}_\varepsilon)} d\bm{x}.
\end{eqnarray}
Since we have used a multivalued cutting line in this paper, getting values for \(\rho\) along the cutting line is actually not possible. Instead, consider the boundary of the fattened cutting line \(\partial \mathcal{D}_\varepsilon\) as \(\varepsilon \to 0\), which is essentially the two sides of the cutting line, which we call \(\mathcal{D}_a\) and \(\mathcal{D}_b\) (the order does not matter here, only that \(a\) and \(b\) are different sides). The integral should then be split into a sided cutting line (each side of width \(\varepsilon\)), such that
\begin{align}
\mathcal{L}^2(\bar{E}) &\approx \lim_{\varepsilon \to 0} \int_{\mathcal{D}} \frac{2\varepsilon}{\rho(\bm{x}, \mathcal{D}_\varepsilon)} d\bm{x} \\
&\approx \lim_{\varepsilon \to 0} \int_{\mathcal{D}_a} \frac{\varepsilon}{\rho(\bm{x}, \mathcal{D}_\varepsilon)} d\bm{x} + \lim_{\varepsilon \to 0} \int_{\mathcal{D}_b} \frac{\varepsilon}{\rho(\bm{x}, \mathcal{D}_\varepsilon)} d\bm{x}.
\end{align}
Since the parameterization in Sec.~\ref{sec:paramet} accounts for the sidedness of the cutting lines already, this can be written succinctly as
\begin{equation}
\mathcal{L}^2(\bar{E}) = \lim_{\varepsilon \to 0} \int_0^{4\pi} \frac{\varepsilon}{\rho(\bm{s}(\theta),\mathcal{D}_\varepsilon)} d\theta. \label{eq:measure_of_E}
\end{equation}

If the regions of intersection of the cutting lines and a particular orbit are known, then the measure of the particular orbit can be computed easily,
\begin{align}
\mathcal{L}^2(\bar{X}) &\approx \lim_{\varepsilon \to 0} \int_{\mathcal{D}_a \cap \bar{X}} \frac{\varepsilon}{\rho(\bm{x}, \mathcal{D}_\varepsilon)} d\bm{x} \nonumber \\
 &\quad + \lim_{\varepsilon \to 0} \int_{\mathcal{D}_b \cap \bar{X}} \frac{\varepsilon}{\rho(\bm{x}, \mathcal{D}_\varepsilon)} d\bm{x}\\
	&\approx [\mathcal{L}^1(\mathcal{D}_a \cap \bar{X}) + \mathcal{L}^1(\mathcal{D}_b \cap \bar{X})] \lim_{\varepsilon \to 0} \frac{\varepsilon}{\rho(\bm{x},\mathcal{D}_\varepsilon)}.
\end{align}
For many orbits, this is zero, but any orbit that is dense in a two-dimensional set will evaluate as the measure of the two-dimensional set.

\begin{figure*}
  \includegraphics[width=0.9\textwidth]{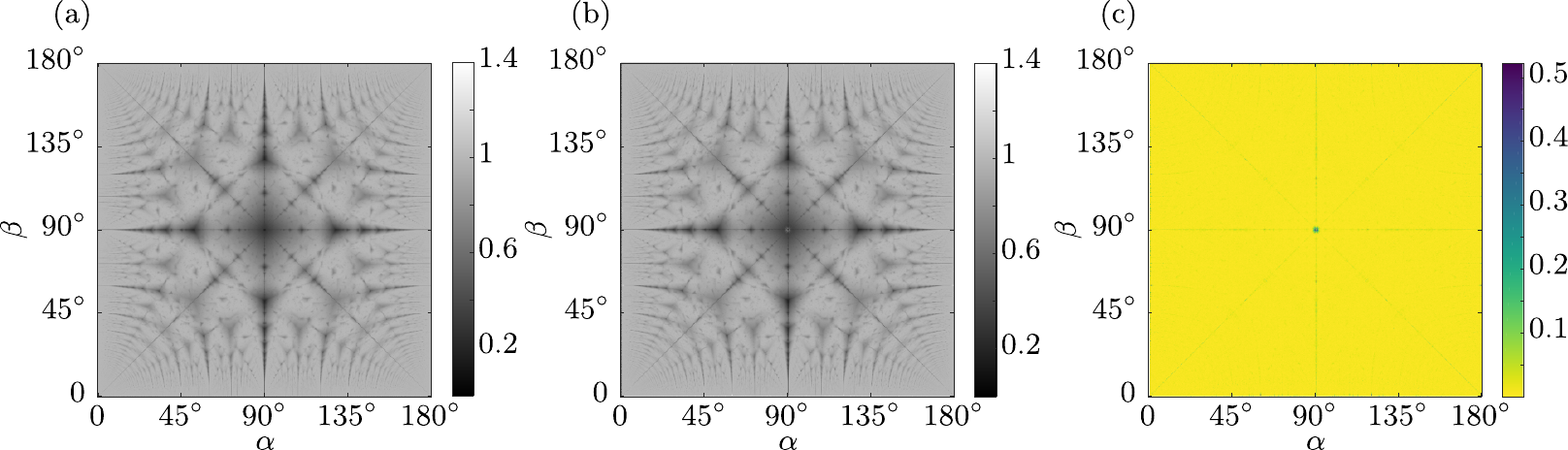}
  \caption{A comparison of \(\Phi\) values computed by (a) direct measurement of the exceptional set as in previous work \cite{Lynn2019} and (b) by integrating the cutting line density. The absolute value of the difference between the two is shown in (c). For most protocols, the difference is near zero. Also note that \(\Phi > 1 \), which is impossible, for some protocols near the boundaries in (b) while values above 1 are impossible using the method in (a). Both methods use \(\varepsilon = 1\times 10^{-3}\) wide cutting lines for \(N = 20\,000\) iterations.
  }
  \label{fig:phase_compare}
\end{figure*}

To illustrate application of this method, the fractional coverage of the hemispherical shell, \(\Phi = \mathcal{L}^2(\bar{E})/2\pi\), was computed using 20\,000 iterations of protocols within the \((\alpha, \beta) \in [0,180^\circ] \times  [0,180^\circ]\) space using the method presented here [Fig.~\ref{fig:phase_compare}(b)] as well a previously used method [Fig.~\ref{fig:phase_compare}(a)] that directly measures a fattened exceptional set \cite{Lynn2019}. Both use \(\varepsilon = 1 \times 10^{-3}\) fattened cutting lines. The direct measurement method samples the fattened cutting line at 1024 evenly spaced points across the shell. The cutting line density measurement uses \(\delta = 1 \times 10^{-6}\) and measures \(\rho(\bm{x}, \mathcal{D}_\varepsilon)\) at 2000 evenly spaced points along the cutting lines. For most protocols, the difference between the two, shown in Fig.~\ref{fig:phase_compare}(c), is near zero, but has large deviations in some regions from the direct measurement method. The large difference near \((90^\circ, 90^\circ)\) could result from orbits being incomplete after 20\,000 iterations. The direct measurement method is less sensitive to incomplete orbits as area accumulates nearly uniformly as iterations increase. The density measurement requires that many circuits be taken around an orbit in order for the value of \(\rho\) to converge near its true value. For most protocols, this happens quickly (i.e.\ within 1000 iterations), but near \((90^\circ, 90^\circ)\), convergence is very slow, which overestimates \(\rho\) and results in an incorrect estimate of \(\Phi\). A more concerning phenomenon occurs near the boundaries when approximating \(\Phi\) using cutting line density, when cutting lines begin to significantly overlap with one another (overlap is roughly \(2\varepsilon^2[\csc(\alpha)\sec(\alpha) + \csc(\beta)\sec(\beta)]\)), which results in an estimate for \(\Phi\) that is larger than 1. The only remedy for this is to use a smaller value of \(\varepsilon\) and a larger number of iterates. Since cutting lines are so close, it becomes challenging to separate the two cutting lines from one another and the estimate breaks down.

The main diagonals in Fig.~\ref{fig:phase_compare} contain protocols that have additional symmetries \cite{Smith2017} and create cells that scale more gradually than protocols off the diagonal. For these protocols, the estimate of \(\Phi\) using the cutting line density is more accurate than the direct measurement, which handles smaller cells sizes poorly, of which there are many along the diagonal.

\section{Analytic density example} \label{sec:app_C_example_density}

For some trivial protocols, the density of cutting lines can be analytically defined. Take the example of the \((\phi,0)\) protocol where \(\phi / \pi \notin \mathbb{Q}\) (i.e.\ irrational). This is an example of irrational rotation that creates an exceptional set dense in \(S\) such that there are no cells, i.e.\ \(O = \emptyset\) where \(O\) is the set of cells as defined in Sec.~\ref{sec:introduction}, and every orbit is dense in an arc across the hemisphere in the \(x\) direction (rotation about the \(z\) axis does not change \(z\) values) with arc-length \(l(z) = \pi \sqrt{1 - z^2}\). For this example, we ignore the cutting line at the equator (due to \(\beta =0\), rotation about the second axis) since it generates no cutting lines. Each orbit is also ergodic such that the density \(\rho(\bm{x},\mathcal{D}_\varepsilon)\) for any point \(\bm{x} = (x,y,z) \in S\), given cutting line width \(\varepsilon > 0\) (on both sides), is the fraction of \(l(z)\) that intersects \(\mathcal{D}_\varepsilon\), which is roughly \footnote{There are some caveats with regard to the curvature of orbit arc versus the curvature of \(\mathcal{D}\) (e.g. the true length is \(2\arcsin(\varepsilon)\)), but these disappear as \(\varepsilon \to 0\).} of length \(2\varepsilon\) due to contributions from both separate sides of the cutting line. Therefore, using sided cutting lines in this analytic example, the cutting line density at any point \(\bm{x} \in S\) is
\begin{equation}
\rho(\bm{x},\mathcal{D}_\varepsilon) = \left\{
	\begin{matrix}
	\frac{2\varepsilon}{l(z)} & \text{ if } 2\varepsilon < l(z) \\
	1 & \text{ if } 2\varepsilon \ge l(z)
	\end{matrix}
	\right..
\end{equation}
Note that in the limit of \(\varepsilon \to 0\), \(\rho \to 0\) everywhere except at the poles where \(\rho\) is always 1. Here, it is immediately clear that \(\rho \sim \varepsilon\) because \(\rho\) can easily be thought of as a fraction of an ergodic orbit, and we have used this in the explicit definition of \(\rho\). The normalized density, i.e.\ \(\rho/\varepsilon\), is simply \(2/l(z)\). In this example, the arc-length of any orbit is closely connected to the normalized density. Density, although symmetric about \(z = 0\), is different for each orbit indicating that all of the domain is confined along infinitely thin orbit curves that do not mix. Rewriting \(l(z)\) so that it can be parameterized by arc-length along \(\mathcal{D}\) using \(\phi \in [0,2\pi)\) for a sided parameterization of the single cutting line gives,
\begin{equation}
l(\phi) =	\pi |\sin(\phi)|,
\end{equation}
and the corresponding fractional coverage using Eq.~\ref{eq:exceptional_set_size_body} normalized by the area of the unit hemisphere,
\begin{eqnarray}
	\Phi  = \frac{\mathcal{L}^2(\bar{E})}{2\pi}  &= & \frac{1}{2\pi} \int_0^{2\pi} \frac{\varepsilon}{\rho}\  d\phi \nonumber \\
	 &= & \frac{1}{2\pi}\int_0^{2\pi} \frac{l(\phi)}{2}\ d\phi \\
	 &= & \frac{1}{2\pi}\int_0^{2\pi} \frac{\pi |\sin(\phi)|}{2} \ d\phi = 1,
\end{eqnarray}
which indicates the entire hemisphere is covered, as expected.

When \(\phi / \pi \in \mathbb{Q}\) (i.e.\ a rational rotation) for the same \((\phi, 0)\) protocol, we expect that \(\Phi = 0\) which is, indeed, the case. First, let \(\phi /\pi = p/q\) such that \(q\) is the periodicity of the entire domain. In this case, there are exactly \(q\) cutting lines of zero width comprising \(\bar{E}\). Then, for small \(\varepsilon > 0\) (with some caveats near the poles that we will ignore since they disappear in the limit of \(\varepsilon \to 0\)), the cutting line density at any point \(\bm{x} \in S\) is
\begin{equation}
	\rho(\bm{x}, \mathcal{D}_\varepsilon) = \left\{
	\begin{matrix}
		1/q, & \text{ if } \bm{x} \in E_\varepsilon \\
		0, & \text{ if } \bm{x} \notin E_\varepsilon
	\end{matrix}\right., 
\end{equation}
where we do not double count multivalued points that exactly return to the cutting line. Each orbit in \(\bar{E}\) can then be assigned the ``arc-length'' \(\varepsilon/\rho\) as before, which is simply \( \varepsilon q\) in this case, and the fractional coverage is
\begin{equation}
	\Phi = \lim_{\varepsilon \to 0} \frac{1}{2\pi} \int_0^{2\pi} \varepsilon q \ d\phi = \lim_{\varepsilon \to 0} \varepsilon q = 0,
\end{equation}
as expected.

The values of \(\Phi = 1\) for irrational rotations and \(\Phi = 0\) for rational rotations are completely unsurprising, but we use this trivial example to demonstrate the concept of density. The closed curves that are the orbits in the irrational case contribute no coverage individually, but together, they cover the hemisphere. For rotation about a single axis, neither the irrational rotation nor the rational rotation cause any mixing since there are no orbits that intersect \(\mathcal{D}\) in more that a finite number of points, and it seems that a condition for mixing in the exceptional set is the existence of orbits that contribute non-zero coverage individually (which would necessarily intersect \(\mathcal{D}\) in some non-zero length region).

\section{Purely IET protocols}\label{sec:app_D_iet_protocols}

\begin{figure}
  \includegraphics[width=0.45\textwidth]{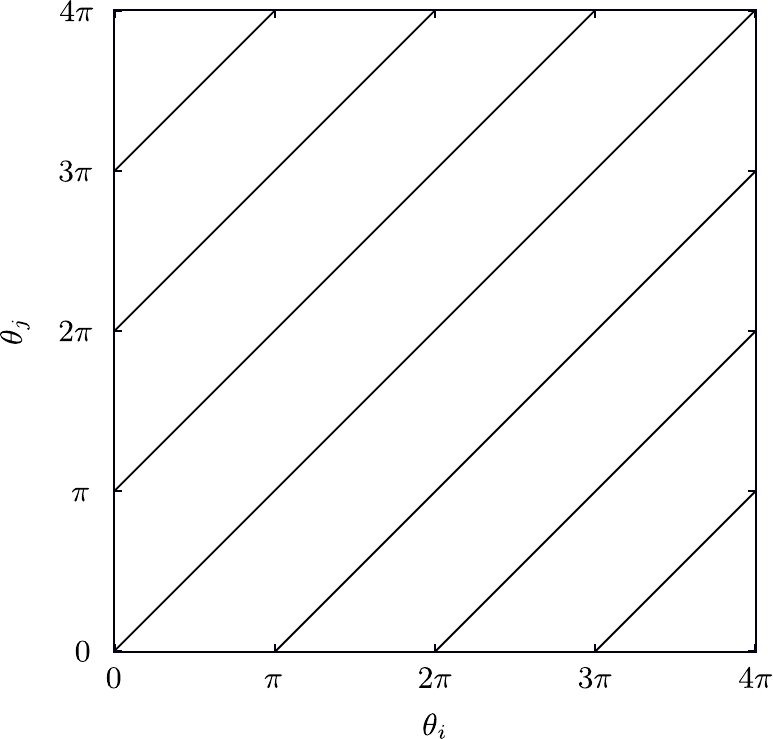}
	\caption{Exact return plot for the \((90^\circ, 90^\circ)\) protocol. This is a protocol with only coincident cutting lines (i.e.\ a polygonal tiling of the sphere) resulting in a periodic interval exchange transform (IET) along the cutting line.}
	\label{fig:90-90-90_return}
\end{figure}

A protocol where \(E\) forms a polygonal tiling (in this case, \(\bar{E}\) is equal to \(E\) since there are no limit points of \(E\)) \cite{Lynn2019}, such as the \((90^\circ, 90^\circ)\) protocol, consists exclusively of coincident (overlapping or repeating) cutting lines and is essentially an interval exchange transform (IET). As such, it is a one-dimensional cut-and-shuffle system \cite{Katok1980, Hmili2010, Keane1977, Masur1982, Avila2007, Krotter2012, Novak2009, Yu2016, Wang2018, Keane1975, Veech1978, Viana2006} on the cutting lines, and the exact return plot can be constructed as shown in Fig.~\ref{fig:90-90-90_return}. For a protocol of this type, the return plot is constructed out of separate diagonal lines consisting of neighboring periodic points along the cutting lines. This results in a continuum of separate invariant sets (which are just collections of periodic points) that makes the entire exceptional set non-mixing and the return plot completely free from mixing blocks.

\end{document}